\magnification=\magstep1
\input amstex
\documentstyle{amsppt}
\catcode`\@=11
\loadmathfont{rsfs}
\def\mycal{\mathfont@\rsfs}
\csname rsfs \endcsname
\catcode`\@=\active

\vsize=6.5in
\topmatter
\title STRONG RIGIDITY OF II$_1$ FACTORS
ARISING FROM MALLEABLE ACTIONS OF $w$-RIGID GROUPS,  I \endtitle
\author SORIN POPA \endauthor

\rightheadtext{Strong rigidity of factors, I}

\affil University of California, Los Angeles\endaffil

\address Math.Dept., UCLA, LA, CA 90095-155505\endaddress
\email popa\@math.ucla.edu\endemail

\thanks Supported in part by a NSF Grant 0100883.\endthanks

\abstract  We consider crossed product II$_1$ factors $M =
N\rtimes_{\sigma} G$, with $G$ discrete ICC groups that contain
infinite normal subgroups with the relative property (T) and
$\sigma$ trace preserving actions of $G$ on finite von Neumann
algebras $N$ that are ``malleable'' and mixing. Examples are the
actions of $G$ by Bernoulli shifts (classical and non-classical),
and by Bogoliubov shifts. We prove a rigidity result for
isomorphisms of such factors, showing the uniqueness, up to
unitary conjugacy, of the position of the group von Neumann
algebra $L(G)$ inside $M$. We use this result to calculate the
fundamental group of $M$, $\mycal F(M)$, in terms of the weights
of the shift $\sigma$, for $G=\Bbb Z^2 \rtimes SL(2, \Bbb Z)$ and
other special arithmetic groups. We deduce that for any subgroup
$S \subset \Bbb R_+^*$ there exist II$_1$ factors $M$ (separable
if $S$ is countable or $S=\Bbb R_+^*$) with $\mycal F(M)=S$.
This brings new light
to a long standing open problem of Murray and von Neumann.
\endabstract
\endtopmatter

\document
\heading 0. Introduction. \endheading

This is the first of a series of papers in which we study rigidity
properties of isomorphisms $\theta$ of {\it crossed product}
II$_1$ factors $M_0, M$ arising from certain actions of groups on
finite von Neumann algebras. We also study isomorphisms between
amplifications of such factors. Typically, we assume the
``source'' factor $M_0$ comes from an action of a group $G_0$
having a large subgroup $H\subset G_0$ with the relative property
(T) of Kazhdan-Margulis ($G_0$ is {\it w-rigid}), while the
``target'' factor $M$ comes from an action $(\sigma,G)$ with good
``deformation+mixing'' properties (a {\it malleable} action), e.g.
an action by Bogoliubov or Bernoulli shifts (classical and
non-classical). The ``ideal'' type of result we seek to prove, is
that any isomorphism between such factors comes from a conjugacy
of the actions involved. Thus, $M_0,M$ can be isomorphic only if
they come from identical group+action data.

Such {\it strong rigidity} results will be obtained in the sequel
papers with the same title [Po6] (the ``group measure space''
case) and [Po7] (the ``non-classical'' case). In the present paper
we prove a key preliminary rigidity result needed in this program,
showing that, after a suitable perturbation, any isomorphism
$\theta$ as above must necessarily take the w-rigid group $G_0$
into the group subalgebra $L(G)$ of $M$. More generally, we prove
that any relatively rigid subalgebra $Q\subset M$ is ``swept'' by
$L(G)$, via a canonical (usually inner) automorphism of $M$.

Besides its r$\hat{\text{\rm o}}$le in ([Po6,7]), this result
already enables us to calculate here the fundamental group $\mycal
F(M)$ of crossed product II$_1$ factors $M$ coming from
(non-classical) Connes-St\o rmer Bernoulli shift actions of
arithmetic groups such as $G=\Bbb Z^2 \rtimes \Gamma$, with
$\Gamma \subset SL(2,\Bbb Z)$ a subgroup of finite index, by using
results from ([Po3], [Ga]). Thus, if $\{t_i\}_i$ are the weights
of $\sigma$ (see [CSt]), then $\mycal F(M)$ is equal to the
multiplicative group generated by the ratios $\{t_i/t_j\}_{i,j}$.
As a consequence, we obtain that any countable subgroup $S \subset
\Bbb R_+^*$ can be realized as a fundamental group of a separable
II$_1$ factor $M$ (i.e., with $L^2(M)$ separable, or equivalently
$M$ countably generated). In fact, by considering Connes St\o rmer
$G$-Bernoulli shifts coming from non-separable Araki-Woods
factors, we obtain II$_1$ factors $M$ with $\mycal F(M)$ any
uncountable subgroup $S\subset \Bbb R_+^*$ as well, but
dim$L^2(M)=|S|$.

This brings new light to a longstanding problem of Murray and von
Neumann on the nature of the fundamental group of II$_1$ factors.
They were led to consider this invariant and to pose this problem
by their theory of continuous dimension and their discovery that
one can take ``$t$ by $t$ matrices'' over a II$_1$ factor $M$, for
any $t>0$. This one parameter family of II$_1$ factors, denoted
$M^t$ and called {\it amplifications} of $M$ by $t$, is used to
define the {\it fundamental group} of $M$ by $\mycal
F(M)=\{t>0\mid M^t \simeq M\}$ ([MvN2]).  After proving that
$\mycal F(R)=\Bbb R_+^*$ (in other words $R^t \simeq R, \forall t
> 0$) for the hyperfinite II$_1$ factor $R$, they comment: ``There
is no reason to believe [that $\mycal F(M)=\Bbb R_+^*$ for all
factors $M$]. The general behavior of this invariant remains an
open question.'' (see [MvN2], page 742). Variants of this problem
were also mentioned in ([K], [Sa], [J3]).

It took almost 40 years until the first progress in this direction
was made, with Connes's breakthrough discovery that for group
factors $M=L(G)$, $\mycal F(L(G))$ reflects the rigidity
properties of the group $G$, being countable whenever $G$ has the
property (T) of Kazhdan. Connes' idea was further exploited in
([Po2], [GoNe], [GeGo], [Po5]) to obtain new classes of separable
II$_1$ factors $M$ with countable $\mycal F(M)$, including
examples for which $\mycal F(M)$ contains a prescribed countable
set. But the first exact computation of a fundamental group $\neq
\Bbb R_+^*$ was obtained only recently, in ([Po3]), where it is
shown, for instance, that $\mycal F(L(G))=\{1\}$ for any group of
the form $G=\Bbb Z^2 \rtimes \Gamma$, with $\Gamma$ a subgroup of
finite index in $SL(2,\Bbb Z)$ (e.g. $\Gamma=\Bbb F_n$). The
result in this paper solves the problem completely, by showing
that any subgroup of $\Bbb R_+^*$ can be realized as a fundamental
group of a II$_1$ factor.

Note however that the problem of finding all subgroups of $\Bbb
R_+^*$ that can occur as fundamental groups of separable II$_1$
factors (which is the really interesting case!) is not completely
solved. Thus, our result only shows that, besides $\Bbb R_+^*$
itself, all countable subgroups of $\Bbb R_+^*$ can appear.
However, we conjecture that the only uncountable subgroup of $\Bbb
R_+^*$ that can occur as a fundamental group of a countably
generated II$_1$ factor is $\Bbb R_+^*$. As a supporting evidence,
note that there are $2^{2^{\aleph_0}}$ many distinct subgroups of
$\Bbb R_+^*$ and only $2^{\aleph_0}$ isomorphism classes of
countably generated II$_1$ factors, so ``most'' uncountable
subgroups of $\Bbb R_+^*$ cannot appear as fundamental groups of
separable factors.

To state in more details the results in this paper, recall some
basic concepts and definitions. Given a finite von Neumann algebra
with a trace $(N,\tau)$ and an action $\sigma: G \rightarrow
\text{\rm Aut}(N,\tau)$ of a discrete group $G$ on $N$ by
$\tau$-preserving automorphisms, its associated {\it crossed
product} von Neumann algebra $N\rtimes_\sigma G$ is generated by a
copy of the group $G$, $\{u_g\}_{g\in G}$, and a copy of the
algebra $N$, acting on the Hilbert space $\Cal H = \oplus_g
L^2(N,\tau) u_g$ by left multiplication subject to the product
rules $u_g \xi u_h=\sigma_{g}(\xi)u_{gh}, \forall g,h \in G$, $\xi
\in L^2(N,\tau)$. Thus, the finite sums $x=\Sigma_g y_g u_g$, $y_g
\in N$ are weakly dense in $N\rtimes_\sigma G$ and in fact any
$\ell^2$-convergent formal sum $x=\Sigma_g y_gu_g$ with
``coefficients'' $y_g$ in $N$ that satisfies $x \cdot \Cal H
\subset \Cal H$ defines an element in $N \rtimes_\sigma G$ (as
left multiplication operator), and all $x\in N\rtimes_\sigma G$
are of this form. The trace $\tau$ on $N$ extends to all $N
\rtimes_\sigma G$ by $\tau(\Sigma_g y_g u_g)=\tau(y_e)$.

This construction goes back to Murray and von Neumann ([MvN1,2]).
The particular case when $N=\Bbb C$ gives the {\it group von
Neumann algebra} $L(G)$ associated to $G$. Thus, $L(G)$ is a
natural subalgebra in any $N \rtimes_\sigma G$. It is a II$_1$
factor iff $G$ is {\it infinite conjugacy class} (ICC). In case
$\sigma$ is an action of $G$ on a probability space $(X,\mu)$ by
measure preserving transformations, it induces an action on the
function algebra $N=L^\infty(X,\mu)$ which preserves $\tau=\int
\cdot {\text{\rm d}}\mu$ and the corresponding crossed product
$L^\infty(X,\mu) \rtimes_\sigma G$ is called the {\it group
measure space} algebra associated to $(\sigma, G)$. It is a II$_1$
factor whenever $\sigma$ is free, ergodic and $G$ infinite.

An example that we often consider in this paper is the (classical)
Bernoulli shift action of $G$ on product spaces $(X,\mu)=\Pi_g
(Y_0,\nu_0)_g$, with base $(Y_0,\nu_0)$ a probability space $\neq$
single point set. These actions are extremely ``malleable'', a
feature that enables us to detect all ``rigid parts'' of  the
group measure space algebra $L^\infty(X,\mu) \rtimes_\sigma G$.
Recall from ([Po3]) that a subalgebra $Q$ of a II$_1$ factor $M$
has the {\it relative property} (T) (or that $Q$ is {\it
relatively rigid} in $M$) if any unital, tracial completely
positive map $\phi$ on $M$ which is close to $id_M$ on a
sufficiently large finite subset of $M$ (in the Hilbert norm
$\|x\|_2 = \tau(x^*x)^{1/2}$) is uniformly close to the identity
on the unit ball of $Q$.

\proclaim{0.1. Theorem} Let $\sigma$ be a Bernoulli shift action
of an ICC group $G$ and denote $M$ the corresponding group measure
space factor. Let $Q\subset M$ be a diffuse von Neumann subalgebra
with the relative property $({\text{\rm T}})$. If $Q$ is either of
type $\text{\rm II}$ or its normalizer in $M$ generates a factor,
then there exists a unitary element $u\in M$ such that $uQu^*
\subset L(G)$. Moreover, if $P$ denotes the von Neumann algebra
generated by the normalizer of $Q$ in $M$ then $uPu^* \subset
L(G)$.
\endproclaim

If $H \subset G_0$ is an inclusion of groups then $L(H) \subset
L(G_0)$ has the relative property (T) iff the pair $(G_0,H)$ has
the relative property (T) of Kazhdan-Margulis ([Ma]; see also
[dHV]). Thus, if $(\sigma_0, G_0)$, $(\sigma,G)$ are actions of
groups on finite von Neumann algebras $(N_0, \tau_0)$, $(N,\tau)$,
$\theta: M_0 \simeq M$ is an isomorphism of the corresponding
crossed product algebras $M_0=N_0 \rtimes_{\sigma_0} G_0$,
$M=N\rtimes_\sigma G$ and $H \subset G_0$ is a subgroup with the
relative property (T), then $Q = \theta(L(H))\subset M$ has the
relative property (T) in $M$. Thus, if we take the groups $G_0, G$
to be ICC and $G_0$ to be {\it weakly rigid} ({\it w-rigid}), i.e.
to have an infinite normal subgroup with the relative property
(T), then Theorem 0.1 shows that any $\theta: M_0 \simeq M$ can be
perturbed by an inner automorphism so that to take $L(G_0)$ into
$L(G)$ (even onto if $G$ is w-rigid as well).

All we actually need for the proof of the above result is the {\it
malleability} of $\sigma: G \rightarrow {\text{\rm Aut}} (N,
\tau)$ (in the case $N$ is abelian, an even weaker property called
{\it sub malleability}, is sufficient). This property amounts to
the existence of an embedding of $N$ as the core of a von Neumann
algebra with discrete decomposition $(\Cal N, \varphi)$ ([C2,3],
[T1,3]) on which $G$ acts by an extension of $\sigma$, such that there
exists a continuous action $\alpha$ of $\Bbb R$ on $\tilde{\Cal
N}=\Cal N \overline{\otimes} \Cal N$ commuting with the product
action $\tilde{\sigma}_g = \sigma_g \otimes \sigma_g$ and
satisfying $\alpha_1(\Cal N \otimes 1) = 1 \otimes \Cal N$. We
call $\tilde{\sigma}$ a {\it gauged extension} for $\sigma$. It
comes with a countable multiplicative subgroup $S(\tilde{\sigma})
\subset \Bbb R_+^*$, given by the almost periodic spectrum of the
discrete decomposition $(\Cal N, \varphi)$ (i.e., of the modular
group associated with $\varphi$). Since the core of $\Cal N
\rtimes G$ is $M = N \rtimes G$, the group $S(\tilde{\sigma})$ is
contained in $\mycal F(M)$. Even more, for each $\beta \in S
(\tilde{\sigma})$, the inclusion $N \rtimes G \subset \Cal
N\rtimes G$ gives rise to a family of $\beta$-scaling
automorphisms $\text{\rm Aut}_\beta(M; \tilde{\sigma})$, any two
of which differ by an inner automorphism of $M$. A stronger
version of this property, called {\it s-malleability}, requires
the existence of an additional period-2 automorphism (``grading'')
$\beta $ of $\tilde{\Cal N}$ that leaves $\Cal N$ pointwise fixed
and satisfies $\beta\alpha_t = \alpha_{-t}\beta, \forall t$.

A classical Bernoulli shift action $\sigma$ of a group $G$ on the
product space $\Bbb T^G$ is easily seen to be s-malleable (sub
s-malleable for arbitrary base space $(Y_0,\nu_0)$) with
$S(\tilde{\sigma})=\{1\}$ and $\text{\rm Aut}(M; \tilde{\sigma})$
coincides with the inner automorphisms of $M$. A non-classical
Connes-St\o rmer Bernoulli shift action $\sigma$ on the infinite
tensor product algebra $(\Cal N, \varphi)= \overline{\otimes}_g
(M_{k\times k}(\Bbb C), \varphi_0)_g$, with $\varphi_0$ a state of
weights $\{t_i\}_i$, is also s-malleable. In this case
$S(\tilde{\sigma})$ is the multiplicative group generated by the
ratios $\{t_i/t_j\}_{i,j} \subset \Bbb R_+^*$ (cf. [AW], [P]). By
replacing the {\it base} $(M_{k\times k}(\Bbb C), \varphi_0)$ of
the $G$-Bernoulli shift in the above construction with an
arbitrary, possibly non-separable IPTF1 Araki-Woods factor $(\Cal
N_0, \varphi_0)$, we still  get a malleable action. In this case,
$S(\tilde{\sigma})$ is equal to the almost periodic spectrum of
$\varphi_0$, and thus can be taken to be any multiplicative
subgroup $S$ of $\Bbb R_+^*$. Actions by weighted Bogoliubov
shifts ([PSt]) satisfy a similar malleability condition. With
these notations we have:

\proclaim{0.2. Theorem} Let $M_i$ be a type ${\text{\rm II}}_1$
factor of the form $M_i = N_i \rtimes_{\sigma_i} G_i$, where $G_i$
is w-rigid ICC, $\sigma_i: G_i \rightarrow {\text{\rm Aut}}(N_i,
\tau_i)$ is a Connes-St\o rmer Bernoulli shift, $i=0,1$. Assume
there exists an isomorphism $\theta : M_0 \simeq M_1^{s}$, for
some $s > 0$. Then there exist $\beta_i \in S(\tilde{\sigma}_i)$
and $\theta^i_{\beta_i} \in \text{\rm Aut}_{\beta_i}(M_i;
\tilde{\sigma}_i)$ such that $\theta^1_{\beta_1}(\theta
(L(G_0)))=L(G_1)^{s\beta_1}$,
$\theta(\theta^0_{\beta_0}(L(G_0)))=L(G_1)^{s\beta_0}$. Moreover,
$\beta_0=\beta_1$ and $\theta^i_{\beta_i}$ are unique modulo
perturbation by an inner automorphism implemented by a unitary of
$L(G_i)$.
\endproclaim

When applied to the case $G_0=G_1=G$, $\sigma_0 = \sigma_1=\sigma$
and $M_0=M_1 = N \rtimes_\sigma G$, with $G$ w-rigid ICC and
$\sigma$ a Connes-St\o rmer Bernoulli shift, the above theorem
shows that if $s \in \mycal F(N \rtimes_\sigma G)$ then there
exists $\beta \in S(\tilde{\sigma})$ such that $s\beta \in \mycal
F(L(G))$. Hence, if $\mycal F(L(G))=\{1\}$ then $s = \beta^{-1}
\in S(\tilde{\sigma})$, implying that $\mycal F(N \rtimes_\sigma
G)=S(\tilde{\sigma})$. Thus, if we take $G = \Bbb Z^2 \rtimes
SL(2,\Bbb Z)$ (which is w-rigid), then the calculation $\mycal
F(L(G))=\{1\}$ in ([Po3]) implies $\mycal F(N\rtimes_\sigma
G)=S(\tilde{\sigma})$. More generally, from the results in ([Po3])
on HT factors and their $\ell^2_{_{HT}}$-Betti numbers
$\beta_n^{^{HT}}$, for which the fundamental group is trivial
whenever there exists some $n$ with $\beta_n^{^{HT}}\neq 0,
\infty,$ one gets :

\proclaim{0.3. Corollary} $1^\circ$. Let $G$ be a w-rigid ICC
group and $\sigma$ a malleable mixing action of $G$ on $(N,
\tau)$, with an extension $\tilde{\sigma}$ having spectrum
$S(\sigma)$. If $\mycal F(L(G))=\{1\}$ then $\mycal
F(N\rtimes_\sigma G)=S(\tilde{\sigma})$. In particular, this is
the case if $L(G)$ is a ${\text{\rm HT}}$ factor and
$\beta^{^{HT}}_n(L(G))\neq 0, \infty$ for some $n$.

$2^\circ$. Let $G_i$ be w-rigid ICC groups and $\sigma_i$
Connes-St\o rmer Bernoulli shifts of $G_i$ on $(N_i, \tau_i)$,
$i=0,1$. Assume $L(G_i)$ are ${\text{\rm HT}}$ factors with
$\beta^{^{HT}}_n(L(G_0))=0$ and $\beta^{^{HT}}_n(L(G_1))\neq
0,\infty$ for some $n$. Then $N_0\rtimes_{\sigma_0} G_0$ is not
stably isomorphic to $N_1\rtimes_{\sigma_1} G_1$.
\endproclaim

It is easy to see that for any countable subgroup $S \subset \Bbb
R_+^*$ there exists a Connes-St\o rmer Bernoulli shift with the
ratios $\{t_i/t_j\}_{i,j}$ of its weights $\{t_i\}_i$ generating
$S$. Thus, any countable subgroup $S\subset \Bbb R_+^*$ can be
realized as $\mycal F(N \rtimes_\sigma G)$, with $G=\Bbb Z^2
\rtimes SL(2,\Bbb Z)$ and $\sigma$ a Connes-St\o rmer Bernoulli
shift action of $G$. Noticing that such actions leave invariant a
Cartan subalgebra of the hyperfinite II$_1$ factor, it follows
that the factor $N \rtimes_\sigma G$ is naturally isomorphic to
the von Neumann algebra of an equivalence relation $\Cal S$ with
fundamental group $\mycal F(\Cal S)\overset\text{\rm def} \to =
\{t > 0 \mid \Cal S^t \simeq \Cal S \}$ equal to $\mycal
F(N\rtimes_\sigma G)$. Moreover, one can take a Connes-St\o rmer
$G$-Bernoulli shift action $\sigma$ on $\overline{\otimes}_g (\Cal
M_0, \varphi_0)_g$, with the base $(\Cal N_0, \phi_0)$ a suitable
non-separable ITPF1 factor, to get $S(\tilde{\sigma})$ to be an
arbitrary uncountable subgroup $S\subset \Bbb R_+^*$ as well. Thus we
have:

\proclaim{0.4. Corollary} Let $S$ be an arbitrary  subgroup of
$\Bbb R_+^*$ and let $\Gamma \subset SL(2,\Bbb Z)$ be a subgroup
of finite index (e.g. $\Gamma \simeq \Bbb F_n$).

$1^\circ$. There exist properly outer actions $\sigma$ of $\Gamma$
on an approximately finite dimensional $\text{\rm II}_1$ factor
$N$ such that $\mycal F(N \rtimes_{\sigma} \Gamma) = S$. Moreover,
if $S$ is countable or equal to $\Bbb R_+^*$ then $N$ can be taken
to be the hyperfinite ${\text{\rm II}}_1$ factor $R$ $($i.e., the
unique approximately finite dimensional ${\text{\rm II}}_1$
factor$)$.

$2^\circ$. If $S$ is countable, 
then there exist countable, measurable, 
measure preserving ergodic standard equivalence relations $\Cal S$ of
the form $\Cal R \rtimes \Gamma$, with $\Cal R$ ergodic
hyperfinite equivalence relation and $\Gamma$ acting outerly on
$\Cal R$, such that $\mycal F(\Cal S) = S$.
\endproclaim

Note that if a II$_1$ factor $M$ comes from an equivalence
relation $\Cal S$ then $\mycal F(\Cal S) \subset \mycal F(M)$.
Related to the problem of showing ``$\mycal F(M)\neq \Bbb R_+^*
\Rightarrow \mycal F(M)$ countable'' for all separable II$_1$
factors $M$, mentioned earlier, it would of course be equally
interesting to prove the similar fact for equivalence relations.

The ideas and techniques used in the proof of 0.1 and 0.2 are
inspired from ([Po1,3,4]). Thus, the malleability of $\sigma$
combined with the the w-rigidity of $G$ allows  a
``deformation/rigidity'' argument in the algebra $\tilde{\Cal M} =
\tilde{\Cal N} \rtimes_{\tilde{\sigma}} G$. As a result of this
argument, we obtain a non-trivial $L(G_0)-L(G)$ Hilbert submodule
of $L^2(\Cal M, \varphi)$ which is finite dimensional as a right
$L(G)$-module ($\Cal M$ denotes here the von Neumann algebra $\Cal
N \rtimes_\sigma G$). Using ``intertwining subalgebras''
techniques similar to ([Po3], A.1), from such a bimodule we get a
unitary element $u \in \tilde{\Cal M}$ that normalizes $M=N
\rtimes G$ and conjugates a corner of $L(G_0)$ onto a corner of
$L(G)$. The trace scaling automorphism $\theta_\beta$ in Theorem
0.2 is then nothing but Ad$u$.

In the case of 0.2, this argument is carried out in the framework
of von Neumann algebras with discrete decomposition (thus possibly
of type III), whose theory was developed in the early 70's by A.
Connes ([C2,3]). Our work benefits directly or indirectly from
these papers, the Tomita-Takesaki theory, Takesaki duality and the
work on type III factors in ([CT], [T1]). However, since our
states $\varphi$ are almost periodic (thus ``almost-like traces'')
the formalism simplifies, allowing us to complete most of the
proofs by just using ``II$_1$-corners'' of II$_\infty$ factors,
their trace scaling automorphisms and the associated crossed
product algebras.

The paper is organized as follows: In Section 1 we recall some
basic facts about discrete decomposition of von Neumann algebras,
introduce definitions and notations related to malleability of
actions, and give examples. In Section 2 we prove several
equivalent conditions for subalgebras of the core $M=\Cal
M_\varphi$ of a factor with discrete decomposition $(\Cal M,
\varphi)$ to be conjugate via partial isometries of $\Cal M$ that
normalize $M$. The effectiveness of these conditions depends on a
good handling of relative commutants for subalgebras in $\Cal
M=\Cal N \rtimes G$, and Section 3 proves the necessary such
results. In Section 4 we prove the main technical result of the
paper: a generalized version of 0.1 showing that if $\sigma$ is
malleable mixing then $L(G)$ ``absorbs'' all relatively rigid
subalgebras of $M= N \rtimes_\sigma G$ (see 4.1 and 4.4).

In Section 5 we prove 0.2-0.4, while in Section 6 we relate the
class of factors studied in this paper with the HT factors
introduced and studied in ([Po3]), showing the two classes are
essentially disjoint. As an application, we prove that
(classical) Bernoulli $\Bbb F_n$-actions
cannot be orbit equivalent to the actions of $\Bbb F_n$
considered in ([Po3]),
thus providing two new free ergodic actions of $\Bbb F_n$,
non-orbit equivalent to the three ones constructed in ([Po3]) and
to the one in ([Hj]).

This work initiated while I was visiting the Laboratoire
d'Alg\'ebres d'Op\'erateurs at the Universit\'e de Paris 7, during
the Summer of 2001. I am grateful to the CNRS and the members of
the Lab for their support and hospitality. During the final stages
of writing this paper I benefitted from many useful comments and
discussions with A. Connes, D. Gaboriau, N. Monod, M. R\o rdam, G.
Skandalis and M. Takesaki, and I take this opportunity to thank
them. The results presented here have been first announced at the
WABASH Conference in September 2002.

\heading 1. Malleable mixing actions and their gauged extensions
\endheading

We first  recall the definition and
properties of von Neumann algebras with
discrete decomposition, which play an important role
in this paper. Then we consider actions of groups on
discrete decompositions and their associated cross product algebras.
Following ([Po1]), we next
define the notion
of gauged extension
and the malleability property for actions. We introduce notations,
list basic properties and give examples.

\vskip .05in \noindent {\bf 1.1. von Neumann algebras with
discrete decomposition.} A von Neumann algebra with a normal
faithful state $(\Cal N, \varphi)$ has a {\it discrete
decomposition} if the centralizer of the sate $\varphi$, $N=\Cal
N_\varphi \overset\text{\rm def} \to = \{x\in \Cal N \mid
\varphi(xy) = \varphi(yx), \forall y \in \Cal N \}$, satisfies
$N'\cap \Cal N = \Cal Z(N)$ and the modular automorphism group
$\sigma^\varphi$ associated with $\varphi$ is almost periodic. The
almost periodic spectrum of $\sigma^\varphi$ is called the {\it
almost periodic spectrum} of the discrete decomposition, and is
denoted $S(\Cal N, \varphi)$. The centralizer algebra $N$ is
called the {\it core} of the discrete decomposition.

The discrete decomposition $(\Cal N, \varphi)$ is
{\it factorial} if $N=\Cal N_\varphi$ is a factor.
It is  {\it tracial} if $N=\Cal N$ and it is
of {\it type} III if $\Cal N$ is a type III von
Neumann algebra.

\vskip .05in
\noindent
{\it 1.1.1. A cross product form for discrete decompositions}.
The structure of a discrete decomposition $(\Cal N, \varphi)$
can be made more specific as follows:

For each $\beta \in \Bbb C$ let $\Cal H^0_\beta = \{x\in \Cal
N\mid \varphi(xy)=\beta \varphi(yx), \forall y\in \Cal N$. Then
$S(\Cal N, \varphi)=\{\beta \mid \Cal H_\beta^0 \neq 0\}$ and the
vector space $\Cal H_\beta^0$ can alternatively be described as
$\{x\in \Cal N \mid \sigma^\varphi_t(x)=\beta^{it}x\}$.

We have $\Cal H^0_\beta \Cal H^0_{\beta'} = \Cal H^0_{\beta\beta'}$ and
$(\Cal H_\beta^0)^* = \Cal H^0_{\beta^{-1}}$,
with $\Cal H^0_1=\Cal N_\varphi=N$ and with
all $\Cal H^0_\beta$ being $N-N$ bimodules. In particular,
$\Sigma_\beta \Cal H^0_\beta$ is a dense $^*$-subalgebra of
$\Cal N$, with $L^2(\Cal H_\beta^0)=\overline{\Cal H^0_\beta}
\subset L^2(\Cal N, \varphi)$,
$\beta \in H(\Cal N, \varphi),$ mutually
orthogonal Hilbert spaces and
$L^2(\Cal N, \varphi)=\oplus_\beta L^2(\Cal H^0_\beta)$.
Thus, if $p_\beta$ denotes the projection onto the
eigenspace $L^2(\Cal H^0_\beta)$, then $\Sigma_\beta p_\beta = 1$.

Note that if $E=E_N$ denotes the $\varphi$-preserving conditional expectation
of $\Cal N$ onto $N$, then $E(\Cal H^0_\beta)=0, \forall \beta \neq 1$.
Also, the projection $p_1$ coincides with the projection
$e_N$ implementing the expectation $E$ by $e_Nxe_N=E(x)e_N$, $x\in \Cal N$
([T2]).

\vskip .05in
\noindent
{\it 1.1.2. The normalizing groupoid}. By the above properties,
if $x \in \Cal H^0_\beta$, for some $\beta \in H(\Cal N, \varphi)$,
and $x=w|x|$ is its polar decomposition, then $|x|\in N$
and the partial isometry $w$ lies in $\Cal H_\beta^0$.
We denote by $\Cal G\Cal V_\beta(\Cal N, \varphi)$ the set
of partial isometries in $\Cal H^0_\beta$ and note that if
$v \in \Cal G\Cal V_\beta(\Cal N, \varphi)$ then
$v^*v, vv^* \in N$, $\tau(vv^*)/\tau(v^*v)=\beta$,
$vNv^* = vv^*Nvv^*$, $\tau(vyv^*)
= \beta \tau(y)$,
$\forall y\in v^*vNv^*v$, where $\tau = \varphi_{|N}$.
Note also that
$\Cal H^0_\beta={\text{\rm sp}}\Cal G\Cal V_\beta(\Cal N, \varphi)$.

\vskip .05in
\noindent
{\it 1.1.3. Fundamental group of the core}.
In case the discrete decomposition $(\Cal N, \varphi)$ is
factorial, each non-zero
$v \in \Cal G\Cal V_\beta (\Cal N, \varphi)$ can be extended to either an
isometry, if $\beta < 1$, or a co-isometry, if $\beta \geq 1$.
Also, if $v \in \Cal G\Cal V_\beta(\Cal N, \varphi)$ then
$\Cal H^0_\beta = {\text{\rm sp}} N v N$.

If $(\Cal N, \varphi)$ is both factorial and infinite dimensional,
then $N$ follows a factor of type II$_1$ and any non-zero $v \in
\Cal G\Cal V_\beta(\Cal N, \varphi)$ implements an isomorphism
$\theta_\beta$ from $N$ onto the $\beta$-amplification $N^\beta$
of $N$, uniquely defined modulo perturbation by inner
automorphisms implemented by unitaries in appropriate reduced
algebras of $N^\infty$. Thus, $S(\Cal N, \varphi)$ is included in
the fundamental group of $N$, $\mycal F(N)$. We denote by
Aut$_\beta(N;\Cal N)$ the set of such isomorphisms $\sigma_\beta$.

\vskip .05in \noindent {\it 1.1.4. Discrete decompositions from
trace scaling actions}. All factorial discrete decompositions
arise as follows (see [C2,3], [T3]): Let $(N, \tau)$ be a type
II$_1$ factor with its unique trace, $Tr=\tau \otimes Tr_{\Cal
B({\ell^2\Bbb N)}}$ the infinite trace on $N^\infty = N
\overline{\otimes} \Cal B(\ell^2\Bbb N)$, $S \subset \mycal F(N)$
a countable subgroup of the fundamental group of $N$ and $\theta$
an action of $S$ on $N^\infty$ by $Tr$-scaling automorphisms,
i.e., such that $Tr \circ \theta_\beta = \beta Tr$, $\forall \beta
\in H$. Let $\Cal N^\infty = N^\infty \rtimes_\theta H$ and
$E^\infty$ the canonical conditional expectation of  $\Cal
N^\infty$ onto $N^\infty$. Let $q=1\otimes q_0$, for some one
dimensional projection $q_0$ in $\Cal B(\ell^2\Bbb N)$. Then
$(\Cal N, \varphi) = (q\Cal N^\infty q, \tau\circ E^\infty_{|\Cal
N})$ has a discrete decomposition and $S = S(\Cal N, \varphi)$,
$E=E^\infty_{|\Cal N}$.

\vskip .05in \noindent {\it 1.1.5. Discrete decomposition for
weights}. The following more general situation will be needed as
well: Let $(\Cal C, \phi)$ be a von Neumann algebra with a normal
semifinite faithful weight. We say that $(\Cal C, \phi)$ has a
{\it discrete decomposition} if $\phi$ is semifinite on the
centralizer von Neumann algebra $\Cal C_\phi = \{ x \in \Cal C\mid
\sigma^\phi_t(x)=x, \forall t\in \Bbb R\}$ and if $(p\Cal C p,
\phi(p \cdot p))$ has discrete decomposition for all projections
$p$ in $\Cal C_\phi$ with $\phi(p) < \infty$. If this is the case,
then we put $S(\Cal C, \phi)\overset \text{\rm def} \to =\cup_p
S(p\Cal C p, \phi(p \cdot p))$.

We denote by $\Cal H^0=\Cal H^0(\Cal C, \phi)$ the Hilbert algebra
$\{x\in \Cal C \mid \phi(x^*x) <\infty, \phi(xx^*) < \infty \}$
and for each $\beta \in S(\Cal C, \phi)$ we let $\Cal
H^0_\beta(\Cal C, \phi) \overset \text{\rm def} \to =\{x\in \Cal
H^0   \mid \phi(xy)=\beta\phi(yx), \forall y\in \Cal H^0\}$. Note
that $x\in \Cal H^0_\beta$ iff $x \in \Cal H^0$ and
$\sigma^\phi_t(x)=\beta^{it} x, \forall t,$ and that  $\Cal H_1^0$
is a hereditary $*$-subalgebra of $\Cal C_\phi$, while
$\Sigma_\beta \Cal H^0_\beta$ is a dense $^*$-subalgebra of $\Cal
C$, with similar multiplicative properties and Hilbert structure
as in the case $\phi$ is a state $\varphi$.

\vskip .05in
\noindent
{\bf 1.2. Actions on discrete decompositions}. Let $(\Cal N, \varphi)$
be a von Neumann algebra with discrete decomposition
and $(N, \tau)=
(\Cal N_\varphi, \varphi_{|\Cal N_\varphi})$
its centralizer algebra, as in Section 1.1.
Let $\sigma : G \rightarrow {\text{\rm Aut}}
(\Cal N, \varphi)$ be a
properly outer action of a discrete group $G$, whose restriction to $(N, \tau)$
is still denoted $\sigma$.

Denote $\Cal M = \Cal N \rtimes_{\sigma} G$ and
$M = N \rtimes_\sigma G$.  We regard $M$ as a subalgebra
of $\Cal M$ in the natural way, with
$\{u_g\}_{g\in G} \subset M \subset \Cal M$ denoting
the canonical unitaries simultaneously
implementing the automorphisms
$\sigma_g$ on $\Cal N, N$.
We denote by $L(G)\subset M \subset \Cal M$ the
von Neumann subalgebra they generate.

We still denote by
$\varphi$ the canonical extension of $\varphi$
from $\Cal N$ to $\Cal M$ and by $E$ the $\varphi$-preserving
conditional expectation of $\Cal M$ onto $M$. Also, $\Cal E$
will denote the canonical ($\varphi$-preserving) conditional
expectation of $\Cal M$ onto $\Cal N$.

Since the action $\sigma$ on $\Cal N$ is $\varphi$- invariant, it
commutes with the corresponding modular automorphism group
$\sigma^\varphi$ on $\Cal N$. Thus, $(\Cal M, \varphi)$ has
discrete decomposition, with core $\Cal M_\varphi=M$. Moreover,
$S(\Cal M, \varphi)=S(\Cal N, \varphi)$, $\Cal H_\beta^0(\Cal
M)=(\Sigma_g \Cal H^0_\beta(\Cal N)u_g)^-= \overline{\text{\rm
sp}} \Cal H^0_\beta(\Cal N)M$ and we have the non-degenerate
commuting square:

$$
\CD
M\ @.\overset{E}\to\subset\ @.\Cal M\\
\noalign{\vskip-6pt}
\cup\ @.\  @.\cup \Cal E\\
\noalign{\vskip-6pt}
N\ @.\subset\ @.\Cal N\\
\endCD
$$

If $\sigma$ is ergodic on the center of $N$ (thus on the center of
$\Cal N$ too, since $\Cal Z(\Cal N) \subset N'\cap \Cal N = \Cal
Z(N)$), then $M=N \rtimes_\sigma G$ follows a factor and, by the
remarks in 1.1, we have $S(\Cal M, \varphi) \subset \mycal F(M)$.

\vskip .05in
\noindent
{\bf 1.3. Basic construction for subalgebras of discrete decompositions}.
We recall here the basic construction associated with
subalgebras of discrete decompositions and explain how
the corresponding extension algebras have discrete decomposition
themselves.
\vskip .05in
\noindent
{\it 1.3.1. The case of subalgebras of the core}.
Let $(\Cal M, \varphi)$ be a von Neumann algebra with  discrete
decomposition and let $B \subset M=\Cal M_\varphi$ be a von Neumann subalgebra.
Let $E_B$ be the unique $\varphi$-preserving conditional
expectation of $\Cal M$ onto $B$,
with $e_B \in \Cal B(L^2(\Cal M, \varphi))$
the orthogonal projection of $L^2(\Cal M, \varphi)$ onto $L^2(B, \varphi)$.
Thus, $e_B$ implements $E_B$ on $\Cal M$ by $e_BXe_B=E_B(X)e_B,
\forall X \in \Cal M$ ([T2]). We denote by
$\langle \Cal M, e_B \rangle$ the von Neumann
algebra generated in $\Cal B(L^2(\Cal M, \varphi))$ by
$\Cal M$ and $e_B$ and call the inclusions
$B \overset{E_B}\to\subset
\Cal M \subset \langle \Cal M, e_B \rangle$
the  {\it basic construction} for $B \subset \Cal M$.

Since $e_B\langle \Cal M, e_B \rangle e_B = Be_B$
and $e_B$ has central support 1 in $\langle \Cal M, e_B \rangle$,
$\langle \Cal M, e_B \rangle$
is an amplification of $B$, thus being semifinite, and
there exists a unique
normal semifinite faithful trace  $Tr$
on $\langle \Cal M, e_B \rangle$
such that $Tr(be_B)=\varphi(b), \forall b\in B$.

However, the ``canonical''
weight on $\langle \Cal M, e_B \rangle$ is not $Tr$
but the following: Let $\Phi$ be the normal
semifinite faithful operator valued weight
of $\langle \Cal M, e_B \rangle$ onto $\Cal M$
determined by $\Phi(xe_B y) = xy,$ $x,y \in \Cal M$, and denote
$\phi = \varphi \circ \Phi$. Then
$\phi$ is a normal semifinite faithful weight on
$\langle \Cal M, \varphi \rangle$ which satisfies
$\phi(e_BY)=\phi(Ye_B)$,
$\phi(xY)=\beta \phi(Yx)$ for all
$x \in \Cal H^0_\beta (\Cal M, \varphi)$,
$Y \in {\text{\rm sp}} \Cal M e_B \Cal M$. Thus, the
modular automorphism group $\sigma^\phi_t$ on $\langle \Cal M, e_B \rangle$
is itself almost periodic and
$(\langle \Cal M, e_B\rangle, \phi)$ has discrete
decomposition with the eigenspace
$\Cal H^0_{\beta_0}(\langle \Cal M, e_B \rangle, \phi)$
being generated by $\Sigma_{\beta}
\Cal H_{\beta_0\beta}^0 e_B \Cal H_{\beta^{-1}}^0$. Moreover,
noticing that $p_\beta \in \langle M, e_B \rangle, \forall \beta$,
$\phi$ is related to the trace $Tr$ by the formula:

$$
\phi(\cdot)=\Sigma_{\beta \in S} \phi(p_\beta \cdot p_\beta)
=\Sigma_{\beta \in S} \beta Tr(p_\beta \cdot p_\beta) \tag 1.3.1
$$

\vskip .05in \noindent {\it 1.3.2. The general case}. Let now
$(\Cal B, \varphi)$ be a discrete decomposition with $\varphi$ a
state and let $\Cal B_0 \subset \Cal B$ be a von Neumann
subalgebra with the property that there exists a
$\varphi$-preserving conditional expectation $E_0$ of $\Cal B$
onto $\Cal B_0$. Let $e_0$ be the orthogonal projection of
$L^2(\Cal B, \varphi)$ onto $L^2(\Cal B_0, \varphi)$ and $\Cal
C=\langle \Cal B, e_0 \rangle$ the von Neumann algebra generated
by $\Cal B$ and $e_0$ on $L^2(\Cal B, \varphi)$. As before,
$\Phi(xe_0y)=xy$ defines a normal semifinite faithful operator
valued weight of $\Cal C$ onto $\Cal B$, with $\phi = \varphi
\circ \Phi$ a normal semifinite faithful weight on $\Cal C$. The
inclusion $\Cal B_0 \subset \Cal B \subset \Cal C=\langle \Cal B,
e_0 \rangle$, with the weight $\phi$ on $\Cal C$, is called the
{\it basic construction} for $\Cal B_0 \overset{E_0}\to\subset
\Cal B$.

We have ${\text{\rm sp}} \Cal C e_0 \Cal C \subset \Cal H^0(\Cal
C, \phi)$ and since $\sigma^\phi_{|\Cal C}=\sigma^\varphi$ and
$\sigma^\phi_t(e_0)=e_0, \forall t$, if $x \in \Cal H^0_\beta(\Cal
B, \varphi)e_0 {\Cal H}^0_{\beta'}(\Cal B, \varphi)$ then
$\sigma^\phi_t(x)=(\beta\beta')^{it} x.$ This shows that the
centralizer Hilbert-algebra $\Cal H_1^0(\Cal C, \phi)$ contains
$\Sigma_\beta \Cal H^0_\beta(\Cal B, \varphi)e_0 {\Cal H}
^0_{\beta^{-1}}(\Cal B, \varphi)$, which has support 1 in $\Cal
C$, and thus $(\Cal C, \phi)$ has discrete decomposition with the
same almost periodic spectrum $S(\Cal C, \phi)$ as $(\Cal B,
\varphi)$.

\vskip .05in \noindent {\bf 1.4. Gauged extensions for actions.}
Let $(N, \tau)$ be a finite von Neumann algebra with a normal,
faithful tracial state. Let $G$ be an infinite discrete group and
$\sigma: G \rightarrow {\text{\rm Aut}} (N, \tau)$ an action of
$G$ on $(N, \tau)$. A {\it gauged extension} for $\sigma$ is an
action $\tilde{\sigma} : G \rightarrow {\text{\rm Aut}}(\Cal N
\subset \tilde{\Cal N}, \tilde{\varphi})$ together with a
continuous action $\alpha: \Bbb R \rightarrow {\text{\rm Aut}}
(\tilde{\Cal N}, \tilde{\varphi})$ satisfying the conditions:

\vskip .05in \noindent $(1.4.1)$. $(\Cal N \subset \tilde{\Cal N},
\tilde{\varphi})$ is an inclusion of von Neumann algebras such
that if we denote $\varphi=\tilde{\varphi}_{|\Cal N}$ and
$\tilde{N}=\tilde{\Cal N}_{\tilde{\varphi}}$ then $N = \Cal
N_{\varphi}=\tilde{N} \cap \Cal N$, $\varphi_{|N}=\tau$, and both
$(\Cal N, \varphi)$, $(\tilde{\Cal N}, \tilde{\varphi})$ have
discrete decomposition, with $S(\Cal N, \varphi)= S(\tilde{\Cal
N}, \tilde{\varphi})$, $\Cal G\Cal V(\Cal N, \varphi) \subset \Cal
G\Cal V(\tilde{\Cal N}, \tilde{\varphi})$ and $\tilde{\Cal N} =
\overline{\text{\rm sp}} \Cal G\Cal V(\Cal N, \varphi)\tilde{N}$.
\vskip .05in \noindent $(1.4.2)$. $\tilde{\sigma}: G \rightarrow
{\text{\rm Aut}} (\tilde{\Cal N}, \tilde{\varphi})$ is an action
such that $\tilde{\sigma}_g(\Cal N)=\Cal N, \forall g,$ and
$\tilde{\sigma}_{| N}= \sigma$. \vskip .05in \noindent $(1.4.3)$.
$\alpha$ commutes with $\tilde{\sigma}$ and satisfies the
conditions:
$$
\overline{\text{\rm sp}}^w \{u\in \Cal U(\alpha_1(\Cal N)) \mid u
\Cal N u^*=\Cal N, {\text{\rm d}}\varphi (u\cdot u^*) /\text{\rm
d}\varphi, {\text{\rm d}}\varphi (u^* \cdot u)/\text{\rm d}\varphi
< \infty \}=\alpha_1(\Cal N) \tag a
$$
$$
\overline{\text{\rm sp}}^w \Cal N \alpha_1(\Cal N)=\tilde{\Cal N}. \tag b
$$
$$
\tilde{\varphi}(y_1\alpha_1(x)y_2) = \tilde{\varphi}(x)
\tilde{\varphi}(y_1y_2), \forall x,y_{1,2}\in \Cal N  \tag c
$$

We denote by $S(\tilde{\sigma})$ the common spectrum $S(\Cal N,
\varphi)=S(\tilde{\Cal N}, \tilde{\varphi})$ and call it the {\it
almost periodic spectrum} of the gauged extension
$\tilde{\sigma}$.

The gauged extension $\tilde{\sigma}$ is {\it tracial} if
$\tilde{\Cal N}=\tilde{N}$ is a finite von Neumann algebra and
$\tilde{\varphi}$ is a trace (afortiori extending the trace $\tau$
of $N$), equivalently if $S(\tilde{\sigma})=\{1\}$.
$\tilde{\sigma}$ is of {\it type} III if $\Cal N, \tilde{\Cal N}$
are of type III and it is {\it factorial} if both $N, \tilde{N}$
are factors.

In Section 4 we will also need gauged extension having an
additional symmetry: A {\it graded gauged extension} for $\sigma$
is a gauged extension $\tilde{\sigma}:G \rightarrow \text{\rm
Aut}(\Cal N \subset \tilde{\Cal N}, \tilde{\varphi})$, $\alpha:
\Bbb R \rightarrow \text{\rm Aut}(\tilde{\Cal N},
\tilde{\varphi})$ together with a period $2$-automorphism
$\beta\in \text{\rm Aut}(\tilde{\Cal N}, \tilde{\varphi})$
satisfying:

$$
\Cal N \subset \tilde{\Cal N}^\beta, \beta \alpha_t =
\alpha_{-t}\beta, \forall t. \tag 1.4.4
$$

Note that if $\sigma_i: G_i \rightarrow {\text{\rm Aut}} (N_i,
\tau_i)$ has gauged extension $(\tilde{\sigma_i}, \alpha_i)$,
$\forall i \in I$, then $\otimes_i \tilde{\sigma}_i$ is a gauged
extension for the action $\otimes_i \sigma_i$ of $\times_i G_i$ on
$(\overline{\otimes}_i N_i, \otimes \tau_i)$, with gauge
$(\otimes_i \alpha_i)_t= \otimes_i \alpha_i(t)$. Moreover,
$S(\otimes_i \tilde{\sigma}_i)$ is the multiplicative subgroup of
$\Bbb R_+^*$ generated by $S(\tilde{\sigma}_i)$, $i\in I$.

Also, if $G_0 \subset G$ is a subgroup and $\tilde{\sigma}$ is a
gauged extension for $\sigma: G \rightarrow {\text{\rm Aut}} (N,
\tau)$, then $\tilde{\sigma}_{|G_0}$ is a gauged extension for
$\sigma_{|G_0}$, with $S(\tilde{\sigma}_{|G_0})
=S(\tilde{\sigma})$.

\vskip .05in \noindent {\bf 1.5. Malleability and mixing
conditions for actions.} An action $\sigma$ of a group $G$ on a
finite von Neumann algebra $(N, \tau)$ is {\it malleable} (resp.
{\it s-mallebale}) if it has a gauged extension (resp. a graded
gauged extension). The action $\sigma$ is {\it malleable mixing}
(resp. {\it s-malleable mixing}) if it is mixing and has a
(graded) gauged extension $\tilde{\sigma}: G \rightarrow
{\text{\rm Aut}} (\Cal N \subset \tilde{\Cal N}, \tilde{\varphi})$
with $\tilde{\sigma}$ mixing, i.e.,

$$
\underset g \rightarrow \infty \to \lim
\tilde{\varphi}(x\tilde{\sigma}_g(y)) = \tilde{\varphi}(x)
\tilde{\varphi}(y), \forall x, y \in \tilde{\Cal N}
$$

By the remarks in 1.2, it follows that if $\sigma_i : G
\rightarrow {\text{\rm Aut}} (N_i, \tau_i)$, $i \in I$, are
malleable (resp. malleable mixing) actions of the same group $G$,
then the diagonal product action $\otimes_i \sigma_i : G
\rightarrow {\text{\rm Aut}} (\overline{\otimes}_i N_i, \otimes_i
\tau_i)$, defined by $(\otimes \sigma_i)_g = \otimes_i
\sigma_{i,g}$, $g \in G$, is a malleable (resp. malleable mixing)
action. Also, note that by ([C1]), if one of the $\sigma_i$'s is
properly outer then the diagonal product action $\otimes_i
\sigma_i$ is properly outer.

\vskip .05in \noindent {\bf 1.6. Examples of malleable mixing
actions}. We show here that the commutative Bernoulli shifts with
diffuse base space are s-malleable, while the Connes-St\o rmer
Bernoulli shifts and the Bogoliubov shifts are all malleable
mixing (see also [Po1]).

\vskip .05in \noindent {\it 1.6.1. Commutative Bernoulli shifts}.
Let $(Y_0, \nu_0)$ be a non trivial probability space. Denote $(X,
\mu)=\Pi_{g} (Y_0, \nu_0)_g$ the infinite product probability
space, indexed by the elements of $G$. It is trivial to check that
$\sigma$ is mixing and properly outer. Such an action is called a
{\it commutative} (or {\it classic}) {\it Bernoulli shift}.

Let us show that if the base space $(Y_0,\nu_0)$ is diffuse, i.e.
$(Y_0, \nu_0)=(\Bbb T, \lambda)$, then $\sigma$ has a graded
gauged extension (so in particular it is malleable mixing). Thus,
put $\Cal N=N=L^\infty(X, \mu)$, $\tau = \int \cdot {\text{\rm
d}}\mu$ and still denote by $\sigma$ the action induced by the
above $\sigma$ on $(N, \tau)$. Then define
$\tilde{\sigma}_g=\sigma_g \otimes \sigma_g$. It is clearly an
action of $G$ on $\tilde{\Cal N} = N\overline{\otimes}N$ that
preserves $\tilde{\tau}=\tau \otimes \tau$. We show that it is a
tracial graded gauged extension for $\sigma$ when identifying $N$
with $N \otimes \Bbb C \subset N \overline{\otimes} N$.

To see this, we first construct a continuous action  $\alpha: \Bbb
R \rightarrow {\text{\rm Aut}}(\tilde{N}, \tilde{\tau})$ commuting
with $\tilde{\sigma}$ and satisfying $\alpha(1)(N \otimes \Bbb
C)=\Bbb C \otimes N$. It is in fact sufficient to construct a
continuous action $\alpha_0 : \Bbb R \rightarrow {\text{\rm Aut}}
(A_0 \overline{\otimes} A_0, \tau_0 \otimes \tau_0)$ such that
$\alpha_0(1)(A_0 \otimes \Bbb C)=\Bbb C \otimes A_0$, where $(A_0,
\tau_0)=(L^\infty(\Bbb T, \lambda), \int \cdot {\text{\rm
d}}\lambda)$. Indeed, because then the product action $\alpha(t)=
\otimes_g (\alpha_0(t))_g$ will do.

To construct $\alpha_0$, let $u$ (resp. $v$) be a Haar unitary
generating $A_0 \otimes \Bbb C \simeq L^\infty(\Bbb T,\lambda)$
(resp. $\Bbb C \otimes A_0$). Thus, $u, v$ is a pair of generating
Haar unitaries for $\tilde{A}_0=A_0 \overline{\otimes} A_0$, i.e.,
$\{u^nv^m\}_{n, m \in \Bbb Z}$ is an orthonormal basis for
$L^2(A_0 \overline{\otimes} A_0, \tau_0 \otimes \tau_0) \simeq
L^2(\Bbb T, \lambda) \overline{\otimes} L^2(\Bbb T, \lambda)$. We
want to construct the action $\alpha_0$ so that
$\alpha_0(1)(u)=v$.

Note that given any other pair of generating Haar unitaries $u',
v'$ for $\tilde{A}_0$, the map $u\mapsto u', v\mapsto v'$ extends
to a $\tilde{\tau}_0=\tau_0\otimes \tau_0$-preserving automorphism
of $\tilde{A}_0$. Also, note that $v, uv$ is a pair of generating
Haar unitaries for $\tilde{A}_0$. Thus, in order to get
$\alpha_0$, it is sufficient to find a  continuous action
${\alpha_0}' : \Bbb R \rightarrow {\text{\rm Aut}} (\tilde{A}_0,
\tilde{\tau}_0)$ such that ${\alpha_0}'(1)(v)=uv$.

Let $h \in \tilde{A}_0$ be a self-adjoint element such that
$exp(2\pi i h)= u$. It is easy to see that for each $t$, $u$ and
$exp(2\pi i th)v$ is a pair of Haar unitaries. Denote by
${\alpha_0}'(t)$ the automorphism $u\mapsto u, v\mapsto exp(2\pi i
th)v$. We then clearly have ${\alpha_0}'(t_1) {\alpha_0}'(t_2) =
{\alpha_0}'(t_1+t_2)$, $\forall t_1, t_2 \in \Bbb R$ and
${\alpha_0}'(1)(v)=uv$.

Finally, to construct the grading $\beta$, we let $\beta_0$ act on
$\tilde{A}_0$ by $\beta(v)=v$, $\beta(u)=u^*$. Then clearly
$\beta_0$ leaves $A_0$ pointwise fixed and $\beta_0 \alpha_{0,t}
=\alpha_{0, -t} \beta_0$. Thus, if we take $\beta= \otimes_g
(\beta_0)_g$ then condition $(1.4.4.)$ is satisfied, showing that
$\sigma$ is s-malleable. .

\vskip .05in \noindent {\it 1.6.2. Non-commutative Bernoulli
shifts}. Let now $(\Cal N_0, \varphi_0)$ be a von Neumann algebra
with discrete decomposition and $(\Cal N,
\varphi)={\overline{\underset g \in G \to \otimes}} (\Cal N_0,
\varphi_0)_g$. It is easy to see that $(\Cal N, \varphi)$ has
itself a discrete decomposition, with spectrum $S(\Cal N,
\varphi)=S(\Cal N_0, \varphi_0)$.

If $x = \otimes_g x_g \in \Cal N$ and $h \in G$ then define
$\sigma_h(x) = \otimes_g x'_g$, where $x'_g = x_{h^{-1}g}, \forall
g$. $\sigma$ is then clearly an action of $G$ on $(\Cal N,
\varphi)$, called the $(\Cal N_0, \varphi_0)$-{\it Bernoulli
shift}. Let further $N=\{x \in \Cal N \mid \varphi(xy) =
\varphi(yx), \forall y \in \Cal N\}$ be the centralizer of the
product state $\varphi$. Thus $\sigma_g(N)=N, \forall g$. The
restriction of $\sigma$ to $N$, still denoted $\sigma$, is called
the {\it Connes-St\o rmer}   $(\Cal N_0, \varphi_0)$-Bernoulli
shift action of the group $G$.

Like in $(1.6.1)$ the actions $\sigma$ of $G$ on $(\Cal N,
\varphi)$ and $N=\Cal N_\varphi$ are properly outer and mixing.
Moreover, in case $(\Cal N_0, \varphi_0)$ is an IPTF1 factor the
action on $N$ follows malleable as well, as shown below (see also
[Po1]):

If $(\Cal N_0, \varphi_0)= (M_{k \times k} (\Bbb C), \varphi_0)$
for some $2 \leq k \leq \infty$, and $\varphi_0$ is the faithful
normal state on $M_{k \times k}(\Bbb C)$ of weights $\{t_j\}_j$,
then $\sigma$ is called the {\it Connes-St\o rmer Bernoulli shift
of weights} $\{t_j\}_j$. Note that in this case $S(\Cal N,
\varphi)=S(\Cal N_0, \varphi_0)$ is equal to the multiplicative
subgroup $S=S(\{t_j\}_j)$ of $\Bbb R_+^*$ generated by the ratios
$\{t_i/t_j\}_{i,j}$, called the {\it ratio group} of $\sigma$.
Notice that $S$ is intrinsic to the construction of $\sigma$.

A gauged extension for
$\sigma$ can be obtained as follows:
Let $(\tilde{\Cal N}, \tilde{\varphi}) = (\Cal N,
\varphi)\overline{\otimes}(\Cal N, \varphi)$ and define
$\tilde{\sigma}: G \rightarrow {\text{\rm Aut}}(\tilde{\Cal N},
\tilde{\varphi})$ by $\tilde{\sigma}_g = \sigma_g \otimes
\sigma_g$.
Let $\{e_{ij}\}_{i,j}$ be matrix units for $M_{k\times k}(\Bbb C)$
chosen so that $\varphi_0$ is given by a trace class operator that
can be diagonalized in Alg$\{e_{ii}\}_i$. It is immediate to check
that
$$
\alpha_0(t) = \Sigma_i e_{ii} \otimes e_{ii}
$$
$$
+ \underset i<j \to \Sigma (\cos \pi t/2 (e_{ii} \otimes e_{jj} +
e_{jj} \otimes e_{ii}) + \sin \pi t/2 (e_{ij} \otimes e_{ji} -
e_{ji} \otimes e_{ij})).
$$
defines a continuous action $\alpha_0$ of $\Bbb R$ on $M_{k \times
k}(\Bbb C) \otimes M_{k \times k}(\Bbb C)$ that leaves
$\varphi_0\otimes \varphi_0$ invariant and satisfies
$\alpha_0(1)(M_{k \times k}(\Bbb C) \otimes \Bbb C) = \Bbb C
\otimes M_{k \times k}(\Bbb C)$. Thus, $\alpha_t = \otimes_g
(\alpha_{0}(t))_g$ defines an continuous  action $\alpha$ of $\Bbb
R$ on $(\tilde{\Cal N}, \tilde{\varphi})$ which commutes with
$\tilde{\sigma}$ and satisfies $\alpha (1)(\Cal N\otimes \Bbb
C)=\Bbb C \otimes \Cal N$.

It is straightforward to check that $\tilde{\sigma}, \alpha$
verify the conditions (1.4.1)-(1.4.3). Thus, $(\tilde{\sigma},
\alpha)$ is a gauged extension for $\sigma$. Since
$\tilde{\sigma}$ is also mixing, $\sigma$ follows malleable
mixing, with $S(\tilde{\sigma})=S$. Also, note that $N=\Cal
N_\varphi$ is in this case isomorphic to the hyperfinite II$_1$
factor $R$.

Let now $S \subset \Bbb R_+^*$ be an arbitrary subgroup and for
each $s\in S$ take $\phi_s$ to be the state on $M_{2 \times
2}(\Bbb C)$ of weights $s(1+s)^{-1}, (1+s)^{-1}$ and let $(\Cal
N_0, \phi_0)=\overline{\otimes}_{s\in S} (M_{2\times 2}(\Bbb C),
\phi_s)_s$. Then the Connes-St\o rmer $(\Cal N_0,
\varphi_0)$-Bernoulli shift action $\sigma$ of $G$ is malleable
mixing and $S(\Cal N_0, \varphi_0)=S$. Indeed, this is because
each $(M_{2\times 2}(\Bbb C), \varphi_s)$-Bernoulli shift is
malleable mixing, so the observations in 1.5 apply.

\vskip .05in
\noindent
{\it 1.6.3. Bogoliubov shifts}. Following ([PSt]), we consider the
following construction:  Let $(\Cal H, \pi)$
be a Hilbert space with a representation
of $G$ on it. Assume $\pi$ satisfies
$$
\pi(g) \not\in \Bbb C1 + \Cal K(\Cal H), \forall g\in G, g\neq e \tag a
$$
$$
\underset g
\rightarrow \infty \to \lim \langle \pi(g)\xi, \eta \rangle = 0, \forall
\xi, \eta \in \Cal H \tag b
$$
An example of such a representation is the left regular representation
of $G$ on $\Cal H = \ell^2(G)$.

Let $CAR(\Cal H)$ be the CAR algebra associated with $\Cal H$
and $\sigma=\sigma_\pi : G \rightarrow {\text{\rm Aut}}
(CAR(\Cal H))$ be the action
on the CAR algebra determined by $\pi$, i.e., $\sigma_g(a(\xi))=
a(\pi(g)(\xi)), \xi \in \Cal H, g\in G$.

For a fixed $0 < t < 1$ denote by $\varphi$ the state on $CAR(\Cal
H)$ determined by the constant operator $t 1$ on $\Cal H$ and
define $(\Cal N, \varphi)$ the von Neumann algebra coming from the
GNS representation for $(CAR(\Cal H), \varphi)$. It is easy to see
that $(\Cal N, \varphi)$ has factorial discrete decomposition,
with $S(\Cal N, \varphi)=\lambda^{\Bbb Z}$, where $\lambda =
t/(1-t)$.

Also, since $\pi$ commutes with the constant
operator $t 1$, $\sigma$ leaves $\varphi$ invariant and thus
can be uniquely extended to an action $\sigma: G \rightarrow
{\text{\rm Aut}}(\Cal N, \varphi)$. Condition $(a)$ implies
$\sigma$ is outer on $(\Cal N, \varphi)$ and condition $(b)$ assures that
$\sigma$ is mixing. Also, take $(\tilde{\Cal N}, \tilde{\varphi})$
to be the GNS construction for $(CAR(\Cal H \oplus \Cal H), \varphi)$,
where $\varphi$ is given by the same constant operator $t 1$, but regarded
on $\Cal H \oplus \Cal H$,
then take $\tilde{\sigma} = \sigma_{\pi\oplus \pi}$ and
$\alpha(t)$ be the action on $(CAR(\Cal H \oplus \Cal H)$ implemented by
the unitaries $cos(\pi it)(e_{11}+e_{22}) + sin(\pi it) (e_{12}-e_{21})$
on $\Cal H \oplus \Cal H$, where $e_{11}$ is the projection onto $\Cal H
\oplus 0$, $e_{22}$ the projection onto $0 \oplus \Cal H$ and
$\{e_{ij}\}_{i,j}$ are matrix units
in the commutant of $\pi\oplus \pi$. It is
straightforward to check that $(\tilde{\sigma}, \alpha)$
is a gauged extension for $\sigma_{|N}$, with
$N = \Cal N_\varphi$  isomorphic to the hyperfinite II$_1$ factor.

\vskip .05in
\noindent
{\bf 1.7. Cross-products associated with
gauged extensions}. Let $\sigma : G \rightarrow {\text{\rm Aut}}
(N, \tau)$ be a properly outer action and
$\tilde{\sigma}: G \rightarrow {\text{\rm Aut}}
(\Cal N \subset \tilde{\Cal N}, \tilde{\varphi})$
a gauged extension for $\sigma$ with gauge $\alpha$.
In addition to the cross-product algebras
$\Cal M = \Cal N \rtimes_{\sigma} G$ and
$M = N \rtimes_\sigma G$ considered in 1.2, we let
$\tilde{\Cal M}
= \tilde{\Cal N} \rtimes_{\tilde{\sigma}} G$ and regard
both $M$, $\Cal M$ as subalgebras of
$\tilde{\Cal M}$, with
$\{u_g\}_{g\in G} \subset M \subset \Cal M
\subset \tilde{\Cal M}$
the canonical unitaries and $L(G)\subset M \subset \Cal M
\subset \tilde{\Cal M}$ the
von Neumann subalgebra they generate. Also, we still denote by
$\tilde{\varphi}$ the canonical extension of $\tilde{\varphi}$
from $\tilde{\Cal N}$ to $\tilde{\Cal M}$ (thus
$\varphi=\tilde{\varphi}_{|\Cal M}$).
Note that since $\alpha$ commutes
with $\tilde{\sigma}$, it implements a continuous
action of $\Bbb R$ on $(\tilde{\Cal M}, \tilde{\varphi})$,
still denoted $\alpha$.

By 1.2, both $(\Cal M, \varphi)$ and $(\tilde{\Cal M},
\tilde{\varphi})$ have factorial discrete decomposition, with
cores $ \Cal M_\varphi = M=N \rtimes_\sigma G$ and respectively
$\tilde{\Cal M}_{\tilde{\varphi}}=\tilde{M}= \tilde{N}
\rtimes_{\tilde{\sigma}} G$, and with $S(\tilde{\Cal M},
\tilde{\varphi})= S(\Cal M, \varphi)=S(\tilde{\sigma})$. For each
$\beta \in S(\tilde{\sigma})$ we denote $\text{\rm
Aut}_\beta(M;\tilde{\sigma})$ the set of isomorphisms
$\theta_\beta : M \simeq M^\beta$ implemented by partial
isometries in $\Cal G\Cal V_\beta(\Cal M, \varphi)$.

\vskip .05in
\noindent
{\bf 1.8. Basic constructions associated with gauged extensions}.
By (1.4.3) b), c) it follows that there exists
a $\tilde{\varphi}$-preserving conditional
expectation $\tilde{E}_1$
of $\tilde{\Cal M}$ onto $\Cal M_1=\alpha_1(\Cal M)$,
defined by $\tilde{E}_1(x_0 \alpha_1(x_1)u_g) =
\tilde{\varphi} (x_0)\alpha_1(x_1)u_g$,
$\forall x_{0,1} \in \Cal N, g\in G$. Also, if $\tilde{e}_1$ denotes
the orthogonal projection of $L^2(\tilde{\Cal M}, \tilde{\varphi})$
onto $L^2(\Cal M_1, \tilde{\varphi})\subset L^2(\tilde{\Cal M},
\tilde{\varphi})$  then $\tilde{e}_1X\tilde{e}_1 =
\tilde{E}_1(X)\tilde{e}_1, \forall X
\in \tilde{\Cal M}$ (cf. [T2]).

Note that ${\tilde{E}}_{1|\Cal M}$ coincides with
the unique $\varphi=\tilde{\varphi}_{|\Cal M}$-preserving conditional
expectation $E_1$ of $\Cal M$ onto $L(G)$. Also, if we denote by
$\langle \tilde{\Cal M},  \tilde{e}_1 \rangle$
the von Neumann
algebra generated in $L^2(\tilde{\Cal M}, \tilde{\varphi})$
by $\tilde{\Cal M}$ and $\tilde{e}_1$,
then $\Cal M_1 \overset{\tilde{E}_1}\to\subset
\tilde{\Cal M} \subset \langle \tilde{\Cal M}, \tilde{e}_1 \rangle$
is a  basic construction, in the sense of 1.3.2. Note that the
von Neumann algebra
$\langle \Cal M, \tilde{e}_1 \rangle$ generated by $\Cal M$ and
$\tilde{e}_1$ is naturally isomorphic to $\langle \Cal M, e_1 \rangle$,
where $E_1=E_{L(G)}$, $e_1=e_{L(G)}$ and
$L(G)\overset{E_1}\to\subset \Cal M \subset \langle \Cal M, e_1\rangle$
is a basic construction as in 1.3.1.

We clearly have $\Cal M_1 = \overline{\text{\rm sp}}^w L(G)\Cal N_1
= \overline{\text{\rm sp}}^w \Cal N_1 L(G)$,
$\tilde{\Cal M} = \overline{\text{\rm sp}}^w  \Cal M\Cal N_1 =$
$\overline{\text{\rm sp}}^w  \Cal N_1\Cal M$,
$\langle \tilde{\Cal M}, \tilde{e}_1 \rangle =$ $\overline{\text{\rm sp}}^w
\langle \Cal M, \tilde{e}_1 \rangle \Cal N_1 =$ $\overline{\text{\rm sp}}^w
\Cal N_1 \langle \Cal M, \tilde{e}_1 \rangle$.
Also, if $\tilde{\Phi}$
denotes the normal semifinite faithful operator valued
weight from $\langle \tilde{\Cal M}, \tilde{e}_1 \rangle$
onto $\tilde{\Cal M}$ defined  by
$\tilde{\Phi}(x\tilde{e}_1 y) = xy, x,y\in \tilde{\Cal M}$, and
$\tilde{\phi}=
\tilde{\varphi} \circ \tilde{\Phi}$ then
$\tilde{\Phi}_{|\langle \Cal M, e_1 \rangle}
= \Phi$ and we have the following non-degenerate
commuting squares:

$$
\CD
\Cal M_1\ @.\overset{\tilde{E}_1}\to\subset\ @.\tilde{\Cal M}\
@.\overset{\tilde{\Phi}}\to\subset
\langle\tilde{\Cal M}, \tilde{e}_1\rangle\\
\noalign{\vskip-6pt}
\cup\ @.\ @.\cup\ @.\cup\\
\noalign{\vskip-6pt}
L(G)\ @.\overset{E_1}\to\subset\ @.\Cal M\ @.
\overset{\Phi}\to\subset
\langle \Cal M, e_1 \rangle\\
\endCD
$$

Noticing that Takesaki's criterion ([T2]) is trivially satisfied,
it follows that there exists a unique $\tilde{\phi}$-preserving
conditional expectation $\Cal F$ of $\langle \tilde{\Cal M},
\tilde{e}_1 \rangle$ onto $\langle \Cal M, e_1 \rangle$. The
expectation $\Cal F$ takes $\tilde{\Cal M}$ (resp. $\Cal M_1$)
onto $\Cal M$ (resp. $L(G)$), with $\Cal F_{|\tilde{\Cal M}}$
(resp. $\Cal F_{|\Cal M_1}$) being the unique
$\tilde{\varphi}$-conditional expectation of $\tilde{\Cal M}$
(resp. $\Cal M_1$) onto $\Cal M$ (resp. $L(G)$).

\heading 2. Intertwining subalgebras in factors with discrete
decomposition  \endheading

Let $(\Cal M, \varphi)$ be a von Neumann algebra with  discrete
decomposition. Thus, $M=\Cal M_\varphi$ is a finite von Neumann
algebra with $\tau=\varphi_{|M}$ its faithful trace and $M'\cap
\Cal M = \Cal Z(M)$. Let $B_0 \subset fMf$, $B \subset M$ be
diffuse von Neumann subalgebras, for some projection $f \in \Cal
P(M)$. We establish in this Section necessary and sufficient
conditions for the existence of partial isometries in $\Cal G\Cal
V(\Cal M, \varphi)$ that ``intertwine'' $B_0$ with subalgebras of
$B$. To this end, we use the notations in 1.3. The proofs are
reminiscent of the proofs of (Theorem A.1 in [Po3]) and (Lemmas 4,
5 in [Po4]).

\proclaim{2.1. Theorem} The
following conditions are equivalent:

$1^\circ$. There exists $a\in B_0'\cap f\langle \Cal M, e_B
\rangle f$, with $a \geq 0$, $a \neq 0$ and $\phi(a) < \infty$.

$2^\circ$. There exist $\beta \in H(\Cal M, \varphi)$ and a
non-zero projection $f_0\in B_0'\cap \langle \Cal M, e_B \rangle$
such that $f_0 \leq p_\beta f$ and $Tr(f_0) < \infty$.

$3^\circ$. There exist $\beta \in H(\Cal M, \varphi)$,
a projection $q_0\in B_0$, a
non-zero partial isometry $v\in \Cal G\Cal V_\beta(\Cal M, \varphi)$
and a non-zero $\xi \in q_0L^2(M)vv^*$
such that if we denote $\xi_0 = \xi v$,
then $q_0B_0q_0\xi_0 \subset \overline{\xi_0 B}$.

$4^\circ$. There exist non-zero projections $q\in B_0$, $p\in B$,
an $($unital$)$ isomorphism $\psi$ of $qB_0q$ into $pBp$ and a
non-zero partial isometry $v_0\in \Cal G\Cal V_{\beta}(\Cal M,
\varphi)$, for some $\beta \in H(\Cal M, \varphi)$, such that
$v_0v_0^* \in (qB_0q)'\cap qMq$, $v_0^*v_0 \in \psi(qB_0q)'\cap
pMp$ and $xv_0=v_0\psi(x), \forall x\in qB_0q$.
\endproclaim
\noindent {\it Proof of} $4^\circ \implies 1^\circ$. If $v_0$
satisfies 4$^\circ$, then $qB_0q$ commutes with $v_0e_Bv_0^*$, so
that if $v_1, v_2, ..., v_n$ are partial isometries in $B_0$ with
$v_i^*v_i \leq q$ and $\Sigma_i v_iv_i^* \in \Cal Z(B_0)$ then $a
= \Sigma_i v_i(v_0e_Bv_0^*)v_i^* \in B_0'\cap \langle \Cal M, e_B
\rangle$, while still $\phi(a) < \infty$ (because $v_i \in B_0
\subset M$  are  in the centralizer of $\phi$). \vskip .05in
\noindent {\it Proof of} $1^\circ \implies 2^\circ$. By 1.3 we
have $p_\beta a p_\beta \in B_0'\cap \langle \Cal M, e_B
\rangle_+$, $\forall \beta \in H = H(\Cal M, \varphi)$, and
$\phi(a)=\Sigma_\beta \phi(p_\beta a p_\beta)$. It follows that
there exists $\beta \in H$ such that $p_\beta a p_\beta \neq 0$.
This implies that all spectral projections $e_t$ of $p_\beta a
p_\beta$ corresponding to intervals $[t, \infty)$ for $t > 0$ lie
in $B_0'\cap \langle \Cal M, e_B \rangle_+$ and satisfy $e_t \leq
p_\beta$, $Tr(e_t) = \beta^{-1} \phi(e_t) \leq (t\beta)^{-1}
\phi(a) < \infty$. Thus, any $f_0=e_t \neq 0$ will satisfy
$2^\circ$.

To prove $2^\circ \implies 3^\circ$ we need the following:

\proclaim{2.2. Lemma} Let $L$ be a
finite von Neumann
algebra acting on the Hilbert space $\Cal H$ and assume
its commutant $L'$ in $\Cal B(\Cal H)$ is also
finite. Let
$L_0 \subset L'$ be a von Neumann subalgebra. Then
there exists a projection $q \in L_0$ and $\xi \in \Cal H$
such that $\xi=q\xi \neq 0$ and $qL_0q \xi
\subset \overline{L\xi}$.
\endproclaim
\noindent
{\it Proof}. If $L$ is of type I then $L_0$ and $L_0'\cap L'$ are
both type I. Let $q$ be a maximal
abelian projection in $L_0$ and $q'$ a maximal
abelian projection in $(L_0q)'\cap qL'q$. Then
$qq'L'qq'$ is abelian, implying that the commutant
of $Lqq'$ in
$\Cal B(\Cal H)$ is abelian. Thus, $Lqq'$ is cyclic in $qq'\Cal H$,
i.e., there exists $\xi \in qq'\Cal H$ such that $\overline{L\xi}
= qq'\Cal H$. In particular, $qL_0q\xi \subset \overline{L\xi}$.

If $L$ has a type I summand, then by cutting with the support
projection of that summand we may assume $L$ itself is type I
and the first part applies.

If $L$ is of type II$_1$ then let $a\hat{\in} \Cal Z(L)=\Cal Z(L')$
be the coupling constant between $L$ and $L'$. Thus, $a$ is a positive
unbounded operator affiliated with $\Cal Z(L)=\Cal Z(L')$.
We have to prove that there exist $q\in \Cal P(L_0)$,
$q'\in \Cal P(L_0'\cap L')$ such that $qq'\neq 0$ and
$ctr_{L'}(qq')\leq a^{-1},$ where $ctr_{L'}$ is the central trace on $L'$.
Indeed, for if so then the coupling constant of
$Lqq'$ on $qq'\Cal H$ is $\leq 1$ and $Lqq'$ is
cyclic in $qq'\Cal H$.

If $e_{[0,1]}(a)\neq 0$ then
the statement follows immediately, by taking $q'
= e_{[0,1]}(a)$ and $q=1$. Thus, we may assume $a \geq 1$ and
by cutting with a projection
in $\Cal Z(L') \subset L_0'\cap L'$
we may also assume $a$ bounded.
If $L_0$ (resp. $L_0'\cap L'$) has a type II$_1$ summand,
then by cutting with the support projection of that summand
we may assume $L_0$ (resp. $L_0'\cap L'$) is of type II$_1$
and then we can find projections $q\in L_0$
(resp. $q'\in L_0'\cap L'$) of arbitrary
scalar central trace in $L_0$ (resp. $L_0'\cap L'$). But then the central
trace of $q$ in $L'$ follows equal to that same scalar,
thus $\leq a^{-1}$ when chosen sufficiently small.

If both $L_0, L_0'\cap L'$ are type I then, by cutting each one
of these algebras by an abelian projection (like in the first part),
we may assume both are abelian. This
implies $L_0'\cap L'$ is a maximal
abelian $^*$-subalgebra of the type II$_1$ von
Neumann algebra $L'$. But then $L_0'\cap L'$ has projections
of arbitrary scalar central trace in $L'$, by ([K2]). By choosing
$q'\in \Cal P(L_0'\cap L')$ of central trace $\leq a^{-1}$,
we are done in this case too.
\hfill Q.E.D.
\vskip .05in
\noindent
{\it Proof of} $2^\circ \implies 3^\circ$.
Let $\Cal H_0 = f_0(L^2(\Cal M, \varphi))
\subset p_\beta(L^2(\Cal M, \varphi))$. Since
$f_0$ commutes with $B_0$ and $Tr(f_0) < \infty$,
it follows that $B_0\Cal H_0 B=\Cal H_0$ and $JBJ'\cap \Cal B(\Cal H_0)$ is
a finite von Neumann algebra. By replacing $\Cal H_0$
by $q(\Cal H_0)$ for some appropriate projection $q$ in
$\Cal Z(B_0)$, we may also assume
dim$(\Cal H_0)_{B} < \infty$,
i.e., the central valued coupling constant
of $JBJ$ in $\Cal B(\Cal H_0)$ is uniformly bounded.

By Lemma 2.2 there exists a projection $q_0 \in B_0$ and a
non-zero vector $\xi_0 \in \Cal H_0$ such that $q_0\xi_0 = \xi_0$
and $q_0B_0q_0 \xi_0 \subset \overline{\xi_0 B}$. Since $\Cal H_0
\subset L^2(\Cal H^0_\beta, \varphi)$, it follows that $\xi_0 =
\xi v$ for some $\xi \in L^2(M, \tau)$ and $v \in \Cal G\Cal
V_\beta(\Cal M, \varphi)$. \vskip .05in \noindent {\it Proof of}
$3^\circ \implies 4^\circ$. Let $\xi_0=\xi v$ with $\xi \in
L^2(M)vv^*$  regarded as a square summable operator affiliated
with $M$. Note that $E_B(v^* \xi^*\xi v) \in L^1(B,\tau)_+$ and
that $\xi_0'=\xi (v E_B(v^* \xi^*\xi v)^{-1/2}v^*) v = \xi_0
E_B(\xi_0^*\xi_0)^{-1/2}$ is still in $L^2(M)v$, satisfies $p_0
\overset \text{\rm def} \to = E_B((\xi_0')^*\xi_0')\in \Cal P(B)$
and
$$
q_0B_0q_0 \xi_0' = q_0B_0q_0 \xi_0E_B(\xi_0^*\xi_0)^{-1/2}
\subset L^2(\xi_0 B)E_B(\xi_0^*\xi_0)^{-1/2}
$$
$$
= L^2(\xi_0 E_B(\xi_0^*\xi_0)^{-1/2}B)
= L^2(\xi_0' B).
$$
Thus, by replacing
$\xi_0$ by $\xi_0'$, we may assume $p_0=E_B(\xi_0^*\xi_0)$ is a
projection in $B$. Also, if $\xi = \xi_0v^* \in L^2(M)$
then $v E_B(\xi_0^*\xi_0)v^*
= E_B(\xi^*\xi) \in \Cal P(B)$.

Let $q\in q_0B_0q_0$ be the minimal projection
with the property that $(q_0-q)\xi_0=0$. We denote
$\psi(x)=E_B(\xi_0^*x\xi_0), x\in qB_0q$, and note that $\psi$ is a
unital, normal, faithful,
completely positive map from $qB_0q$ into $p_0Bp_0$.

Also, since $x \xi_0 \in L^2(\xi_0 p_0Bp_0) = \xi_0 L^2(p_0Bp_0)$,
it follows that
$x \xi_0 = \xi_0 \psi(x), \forall x \in qB_0q$. Indeed,
for if $x \xi_0 = \xi_0 y$, for some $y \in L^2(p_0Bp_0)$, then
$\xi_0^*x\xi_0 = \xi_0^*\xi_0 y$ and
so
$$
\psi(x)= E_B(\xi_0^*x\xi_0) = E_B(\xi_0^*\xi_0 y)=
E_B(\xi_0^*\xi_0)y=y.
$$
Thus,
for $x_1, x_2 \in qB_0q$ we get $x_1 x_2 \xi_0 = x_1 \xi_0 \psi(x_2)
= \xi_0 \psi(x_1) \psi(x_2)$. Since we also have
$(x_1x_2) \xi_0 = \xi_0 \psi(x_1x_2)$, this shows that
$\psi(x_1x_2) = \psi(x_1)\psi(x_2)$. Thus, $\psi$ is a unital
$*$-isomorphism of $qB_0q$ into $p_0Bp_0$.

Thus, since $x \xi_0 = \xi_0 \psi(x)$ and $\xi_0^* x = \psi(x)
\xi_0^*, \forall x \in qB_0q,$ it follows that $[qB_0q,
\xi_0\xi_0^*] =0$. Since $\xi_0 = \xi v$, $[qB_0q, \xi \xi^*] =0$
as well. Thus, if $w=(\xi \xi^*)^{-1/2}\xi$ then $w$ is a partial
isometry in $M$ and $x wv = wv \psi(x)$, $\forall x \in qB_0q$.
Thus $v_0 = wv$ will do.
\hfill Q.E.D.

\proclaim{2.3. Corollary} Assume condition $2.1.4^\circ$ is not
satisfied. $($Note that this is the case if there exists no
embedding $\theta: p_0B_0p_0 \hookrightarrow B$, for non-zero $p_0
\in \Cal P(B_0)$.$)$ Then we have:
$$
\forall a_1, a_2, ..., a_n \in M = \Cal M_\varphi, \forall
\varepsilon >0, \exists u \in \Cal U(B_0), \|E_B(a_iua_j^*)\|_2
\leq \varepsilon, \forall i,j. \tag 2.3.1
$$
If in addition $\Cal M$ is finite, with $\varphi$ its trace, then
conversely, $(2.3.1)$ implies non-$2.1.4^\circ$.
\endproclaim
\noindent {\it Proof}. Assume by contradiction that there do exist
$a_1, a_2, ..., a_n \in M$ and $c > 0$ such that $\Sigma_{i,j}
\|E_B(a_iua_j^*)\|^2_\varphi \geq c$, $\forall u\in \Cal U(B_0)$,
and let $b= \Sigma_i a_i^* e_B a_i \in {\text{\rm sp}} M e_B M
\subset \langle \Cal M, e_B \rangle$. We then have the estimates:
$$
Tr(bubu^*)=\Sigma_{i,j} Tr(a_i^*e_Ba_iua_j^*e_Ba_ju^*)
$$
$$
=\Sigma_{i,j} Tr(e_Ba_ju^*a_i^*e_Ba_iua_j^*e_B) =\Sigma_{i,j}
Tr(E_B(a_ju^*a_i^*)E_B(a_iua_j^*))
$$
$$
= \Sigma_{i,j} \|E_B(a_iua_j^*)\|^2_2 \geq c,
$$
for all $u\in \Cal U(B_0)$. Let then $a$ be the element of minimal
norm $\|\cdot \|_{2,Tr}$ in the weak closure of co$\{ubu^* \mid u
\in \Cal U(B_0)\}$ in $\langle \Cal M, e_B \rangle$. Thus, $0\leq
a \leq 1$, $Tr(a) \leq Tr(b)$ and $a \in B_0'\cap
\overline{\text{\rm sp}}^w Me_BM $. Also, $Tr(ba) \geq c > 0$,
implying that $a \neq 0$. By $2^\circ \implies 4^\circ$ in Theorem
2.1, it follows that there exists an isomorphism $\theta:
p_0B_0p_0 \hookrightarrow B$, for some $p_0 \in \Cal P(B_0),
p_0\neq 0$, a contradiction.

To prove the converse, note that if $2.1.2^\circ$ holds true and
$f_0 \in B_0'\cap \angle M, e_B \rangle$ is a finite projection
with $Tr(f_0) < \infty$ then by (1.4 in [P5]) we may assume $f_0 =
\Sigma_j a_j e_B a_j^*$ for some finite set $a_1, a_2, ..., a_n
\in M$, which in turn implies $$ \Sigma_{i,j}
\|E_B(a_iua_j^*)\|^2_2=Tr(f_0uf_0u^*)=Tr(f_0), \forall u \in \Cal
U(B_0),
$$
thus contradicting $(2.3.1)$. \hfill Q.E.D.

\heading 3. Controlling intertwiners and relative commutants
\endheading

Theorem 2.1 shows the importance of controlling intertwiners and
relative commutants of subalgebras of a factor when having to
decide whether the subalgebras are conjugate or not. We prove in
this section two results along these lines:

\proclaim{3.1. Theorem} Let $(\Cal N, \varphi)$ be a von Neumann
algebra with discrete decomposition, $G$ an infinite discrete
group and $\sigma : G \rightarrow {\text{\rm Aut}} (\Cal N,
\varphi)$ a properly outer mixing action. If $Q_0 \subset L(G)$ is
a diffuse von Neumann subalgebra and $x\in \Cal M= \Cal
N\rtimes_\sigma G$ satisfies $Q_0 x \subset \Sigma_i  x_i L(G)$,
for some finite set $x_1, x_2, \ldots, x_n \in \Cal M$, then $x\in
L(G)$. In particular, $Q_0'\cap \Cal M\subset L(G)$ and if $v \in
\Cal M$ is a partial isometry with $[v^*v, Q_0]=0$ and
$vQ_0v^*\subset L(G)$ then $v \in L(G)$.
\endproclaim

\proclaim{3.2. Theorem} Let $\sigma: G \rightarrow {\text{\rm
Aut}}(N, \tau)$ be a malleable mixing action with gauge extension
$\tilde{\sigma}: G \rightarrow {\text{\rm Aut}}(\Cal N \subset
\tilde{\Cal N}, \tilde{\varphi})$. Denote $\tilde{N}=\tilde{\Cal
N}_{\tilde{\varphi}}$, $M= N \rtimes_\sigma G$, $\tilde{M} =
\tilde{N} \rtimes G$, as in $1.8$. Let $P_0 \subset M$ be a
diffuse von Neumann subalgebra such that no corner of $P_0$ can be
embedded (non-unitally) into $N$. Then $P_0'\cap \tilde{M} \subset
M$.
\endproclaim

Both these theorems will be derived from a  general technical
result. To state it we need some notations. Thus, we let $(\Cal T,
\varphi)$ be a von Neumann algebra with discrete decomposition,
$G$ an infinite discrete group, $\sigma : G \rightarrow {\text{\rm
Aut}} (\Cal T, \varphi)$ a properly outer action, $\Cal B = \Cal T
\rtimes_\sigma G$ its cross product algebra with canonical state
$\varphi$, as in $1.2$. Let $\Cal T_0 \subset \Cal T$ be a
$\sigma$-invariant von Neumann subalgebra on which there exists a
$\varphi$-preserving conditional expectation $E_0$. Denote $\Cal
B_0 = \Cal T_0 \rtimes_\sigma G \subset \Cal B$ and still denote
by $E_0$ the $\varphi$-preserving expectation of $\Cal B$ onto
$\Cal B_0$ extending the expectation of $\Cal T$ onto $\Cal T_0$.
Let $\Cal B_0 \overset{E_0}\to\subset \Cal B \subset \Cal
C=\langle \Cal B, e_0 \rangle$ be the basic construction for $\Cal
B_0 \overset{E_0}\to\subset \Cal B$ with its canonical weight
$\phi$, as in $1.3.2$.

\proclaim{3.3. Proposition} Assume there exists $\{b_n\}_n\subset
\Cal T \ominus \Cal T_0$ such that $\{1\} \cup \{b_n \}_n$ is an
orthonormal basis of $\Cal T$ over $\Cal T_0$ with each $b_n$ in
some $\Cal H^0_{\beta_n}(\Cal T, \varphi)$, $\forall n$, and such
that the following condition is satisfied:
$$
\underset g \rightarrow \infty \to \lim ({\text{\rm sup}}\{ \|
E_0(b^*_i y \sigma(g)(b_j))\|_\varphi \mid y \in \Cal N_0, \|y\|
\leq 1 \}) = 0, \forall i,j. \tag 3.3.1
$$
Let $P_0 \subset (\Cal B_0)_\varphi$ be a diffuse von Neumann
subalgebra satisfying the property:
$$
\forall K \subset G \quad finite, \quad \forall \delta > 0,
\exists u\in \Cal U(P_0) \quad with \quad \| \Cal
E(uu_h^*)\|_\varphi \leq \delta, \forall h \in K. \tag 3.3.2
$$
If $\Cal H^0_1=\Cal H^0_1(\Cal C, \phi)$ denotes the centralizer
Hilbert algebra (see $1.3.2$) and $a \in P_0'\cap \Cal H^0_1$ then
$e_0ae_0=a$.
\endproclaim

To prove 3.3, we first need the following:

\proclaim{3.4. Lemma} Under the hypothesis of $3.3$, for any $n$
and any $\varepsilon > 0$ there exists a finite subset $K \subset
G$ and $\delta > 0$ such that if $u \in \Cal U(\Cal B_0)$
satisfies $\|\Cal E(uu_h^*)\|_\varphi \leq \delta$, $\forall h \in
K$, then $\|E_0(b_i^*ub_j)\|_\varphi \leq \varepsilon$, $\forall
i,j$.
\endproclaim
\noindent {\it Proof}. Let $u = \Sigma_g y_g u_g$, with $y_g \in
\Cal T_0$. Then $E_0(b_i^* y_g u_g b_k) = E_0(b^*_i y_g
\sigma_g(b_j)) u_g$ implying that
$$
\|E_0(b_i^*ub_j)\|_\varphi^2 = \Sigma_g \|E_0(b^*_i y_g
\sigma_g(b_j)) \|^2_\varphi.
$$

By $(3.3.1)$ there exists a finite subset $K \subset G$ such that
${\text{\rm sup}}\{ \| E_0(b^*_i y \sigma(g)(b_j))\|_\varphi \mid
y \in \Cal N_0, \|y\| \leq 1 \}\leq \varepsilon/2$, $\forall g \in
G \setminus K$. On the other hand, since the norm $\|\cdot
\|_\varphi$ implements the strong operator topology on the unit
ball of $\Cal M$ and the maps $\Cal N_0 \ni y \mapsto b^*_i y
\sigma(g)(b_j)$ are continuous with respect to the strong operator
topology, $\forall i,j$, it follows that there exists $\delta > 0$
such that if $y \in \Cal N_0$, $\|y\| \leq 1$, $\|y\|_\varphi \leq
\delta$, then $\| E_0(b^*_i y \sigma(g)(b_j))\|_\varphi \leq
(2|K|)^{-1}\varepsilon$. Thus, if $u$ satisfies
$\|y_h\|_\varphi=\|E_0(uu_h^*)\|_\varphi \leq \delta$, $\forall h
\in K$, then

$$
\|E_0(b_i^*ub_j)\|_\varphi^2 = \Sigma_{g\notin K} \|E_0(b^*_i y_g
\sigma_g(b_j)) \|^2_\varphi + \Sigma_{h\in K} \|E_0(b^*_i y_h
\sigma_g(b_j)) \|^2_\varphi
$$
$$
\leq \varepsilon/2+\varepsilon/2 =\varepsilon, \forall i,j.
$$
\hfill Q.E.D.

\noindent {\it Proof of Proposition 3.3}. Recall from 1.3.2 that
$\Cal H^0_1$ is hereditary and contains the $^*$-algebra
$\Sigma_\beta \Cal H^0_\beta(\Cal B, \varphi) e_0 \Cal
H^0_\beta(\Cal B, \varphi)^*$. For $X \in \Cal H^0(\Cal C, \phi)$
(the Hilbert algebra of $(\Cal C, \phi)$, as in 1.1.5 and 1.3.2),
denote $\|X\|_{2, \phi}=\phi(X^*X)^{1/2}$.

Since $\Cal H^0_1=\Cal H^0_1(\Cal C, \phi)$ is a $^*$-algebra and
$(1-e_0)\Cal H^0_1(1-e_0) \subset \Cal H^0_1$, it follows that if
$(1-e_0)a \neq 0$ (resp. $a(1-e_0)\neq 0$), then  by replacing $a$
by a spectral projection of $(1-e_0)aa^*(1-e_0)$ (resp.
$(1-e_0)a^*a(1-e_0)$) corresponding to some interval $[c, 1]$ with
$c > 0$, we may assume $a=f\neq 0$ is a projection with  $f \leq
1-e_0$.

Let $\varepsilon > 0$. Since $\{1\}\cup \{b_n\}_n$ is an
orthonormal basis of $\Cal T$ over $\Cal T_0$, it is also an
orthonormal basis of $\Cal B$ over $\Cal B_0$. Thus, there exists
$n$ such that the orthogonal projection $f_0=\Sigma_{j \leq n} b_j
e_0 b_j^*$ of $L^2(\Cal B, \varphi)$ onto the closure of
$\Sigma_{j \leq n} b_j \Cal B_0$ in $L^2(\Cal B, \varphi)$
satisfies $\|f_0f-f\|_{2, \phi} \leq \varepsilon/3$. Since each
$b_i$ lies in some $\Cal H^0_{\beta_i}(\Cal T,\varphi)$, it
follows that $f_0$ lies in $\Sigma_\beta \Cal H^0_\beta(\Cal B,
\varphi) e_0 \Cal H^0_\beta(\Cal B, \varphi)^*\subset \Cal
H^0_1(\Cal C,\phi)$.

Thus, if $u\in \Cal U(P_0)$ then
$$
\|uf_0u^*f-f\|_{2,\phi} = \|u(f_0f-f)u^*\|_{2,\phi} = \|f_0f-f\|_{2,\phi}
\leq \varepsilon/3.
$$
implying that
$$
\|uf_0u^*f - f_0f\|_{2,\phi} \leq 2 \|f_0f-f\|_{2,\phi} \leq 2\varepsilon/3.
$$

Since $f, f_0$ are in the centralizer Hilbert algebra $\Cal H^0_1$
and $\phi(uXu^*)=\phi(X)$, $\forall X \in \Cal H^0_1, u\in P_0$,
by the Cauchy-Schwartz inequality we get
$$
|\phi(ff_0uf_0u^*f)|= |\phi(ff_0uf_0u^*)|
\leq \|f\|_{2,\phi}  \|f_0uf_0u^*\|_{2,\phi}
$$
$$
= \phi(f)^{1/2} \phi(uf_0u^*f_0uf_0u^*)^{1/2} = \phi(f)^{1/2}
\phi(f_0u^*f_0u)^{1/2}.
$$

Since $\phi(b_iX)=\beta_i\phi(Xb_i)$, $\forall X \in \Cal H^0_1$,
$b_i \in \Cal H^0_{\beta_i}(\Cal B, \varphi)$ and $\|uf_0u^*f -
f_0f\|^2_{2,\phi} = 2\|f_0f\|_{2,\phi}^2 - 2 {\text{\rm Re}} \phi
(ff_0uf_0u^*f)$, by the definition of $\phi$ we obtain
$$
2\|f_0f\|_{2,\phi}^2 \leq \|uf_0u^*f - f_0f\|^2_{2,\phi} +
2\phi(f)^{1/2} \phi(f_0u^*f_0u)^{1/2}
$$
$$
\leq (2\varepsilon/3)^2 + 2\phi(f)^{1/2}\Sigma_{1\leq i,j \leq n}
\varphi(b_iE_0(b_i^*u^*b_j) b_j^*u)
$$
$$
= (2\varepsilon/3)^2 + 2\phi(f)^{1/2}\Sigma_{i,j} \beta_i
\varphi(E_0(b_i^*u^*b_j) b_j^*ub_i)
$$
$$
=(2\varepsilon/3)^2 + 2\phi(f)^{1/2}\Sigma_{i,j} \beta_i
\|E_0(b_j^*ub_i)\|^2_\varphi.
$$

By Lemma 3.4 there exist a finite subset $K \subset G$ and $\delta
> 0$ such that if $u \in \Cal U(P_0)$ satisfies $\|\Cal
E(uu_h^*)\|_\varphi\leq \delta$, $\forall h \in K$, then
$\|E_0(b_i^*ub_j)\|_\varphi \leq
\varepsilon^2/(6n^2\beta\phi(f)^{1/2})$, $\forall 1\leq i,j \leq n
$, where $\beta=\text{\rm max}_i \{\beta_i \mid 1 \leq i \leq
n\}$. By condition $(3.3.2)$ applied for this $K$ and $\delta$,
there exists $u\in \Cal U(P_0)$ such that $\|\Cal
E(uu_h^*)\|_\varphi\leq \delta, \forall h\in K$.

Altogether, since the above summations $\Sigma_{1 \leq i,j \leq
n}$ have $n^2$ terms, we get the estimates
$$
2\|f_0f\|_{2,\phi}^2 \leq (2\varepsilon/3)^2 + 2\phi(f)^{1/2} n^2
\beta \varepsilon^2/(6n^2\beta\phi(f)^{1/2})
$$
$$
\leq (2\varepsilon/3)^2
+ \varepsilon^2/3 \leq
7\varepsilon^2/9.
$$

Thus, $\|f\|_{2,\phi} \leq \|f-ff_0\|_{2,\phi} +\|ff_0\|_{2,\phi}
\leq \varepsilon/3 + 2\varepsilon/3 =\varepsilon$. Since
$\varepsilon > 0$ was arbitrary, this shows that $f=0$, a
contradiction. \hfill Q.E.D.

\vskip .05in \noindent {\it Proof of Theorem 3.1}. Let $L(G)
\overset{E_1}\to\subset \Cal M \subset \langle \Cal M, e_1
\rangle$ be the   basic construction corresponding to the
$\varphi$-preserving conditional expectation $E_1$ of $\Cal M$
onto $L(G)$, like in 1.3.1, with $E_1=E_{L(G)}$ and
$e_1=e_{L(G)}$. We also consider the partition of 1 given by the
projections $\{p_\beta\}_\beta$, as in 1.3.1.

Let $f_x$ be the orthogonal projection of $L^2(\Cal M, \varphi)$
onto the closure $\Cal H$ of $\widehat{Q_0xL(G)}$ in $L^2(\Cal
M,\varphi)$. Since $Q_0 \Cal H L(G)=\Cal H$, $f_x \in Q_0'\cap
\langle \Cal M, e_1 \rangle$. Also, since $\Cal H$ is contained in
the closure of $\widehat{\Sigma_i x_i L(G)}$, it follows that
$\Cal H$ is a finitely generated right $L(G)$-Hilbert module. Thus
$f_x$ lies in the ideal $\Cal J_{e_1}$ generated in $\langle \Cal
M, e_1 \rangle$ by $e_1$. In particular $Tr(f_x) < \infty$, where
$Tr$ is defined as in 1.3.1.

Since $\hat{1}\in L^2(\Cal M, \varphi)$
is a separating vector for $\Cal M$, to prove that $x \in L(G)$
it is sufficient to show that $f_x \leq e_1$.

Assuming $(1-e_1)f_x(1-e_1) \neq 0$, it follows that either
$a=p_\beta f p_\beta\neq 0$, for some $\beta \neq 1$, or
$a=(p_1-e_1)f(p_1-e_1) \neq 0$. Since $p_\beta \in M'\cap \langle
\Cal M, e_1 \rangle$, it follows that $a \in Q_0'\cap \langle \Cal
M, e_1 \rangle$, $a \in \Cal J_{e_1}$, $Tr(a) < \infty$. Also,
since under $p_\beta$ the weight $\phi$ defined by
$\phi(y_1e_1y_2)=\varphi(y_1y_2), y_{1,2} \in \Cal M$ is
proportional to $Tr$ (cf 1.3.1), it follows that $a$ is in the
centralizer of $\phi$. Thus, all spectral projections $f$ of $a$
corresponding to  intervals of the form $[c,\infty)$ for $c > 0$
will satisfy $f\in Q_0'\cap \langle \Cal M, e_1 \rangle$, $f \in
\Cal J_{e_1}$, $\phi(f) < \infty$, $f \leq 1-e_1$ and $f$ in the
centralizer of $\phi$.

Thus, if we apply Proposition 3.3 for $\Cal T=\Cal N$, $\Cal
T_0=\Bbb C$, $\Cal T_1 = \Cal N$ and $\Cal B_0=B_0 = L(G)$, then
we get $f=0$, a contradiction. \hfill Q.E.D. \vskip .05in
\noindent {\it Proof of Theorem 3.2}. Put $\Cal T=\tilde{\Cal N},
\Cal B = \tilde{\Cal M}$, $\Cal T_0 = \Cal N$, $T_0=N$, $\Cal
B_0=\Cal M$, $B_0=M$. Conditions (b) and (c) of $(1.4.3)$ show the
existence of a $\tilde{\varphi}$-preserving conditional
expectation $E_0$ of $\Cal B=\tilde{\Cal M}$ onto $\Cal B_0=\Cal
M$: On elements of the form $y_1 \alpha_1(y_2)$ with $y_1,y_2 \in
\Cal N,$ which by $(1.4.3)$ are total in $\tilde{\Cal N}$, it acts
by $E_0(y_1 \alpha_1(y_2))=\tilde{\varphi}(y_2)y_1$.

By Corollary 2.3, condition $(3.3.2)$ is satisfied. Let us show
that there exists an orthonormal basis $\{1\}\cup \{b_n\}_n$ of
$\Cal T=\tilde{\Cal N}$ over $\Cal T_0 = \Cal N$ verifying
$(3.3.1)$. By the definition of $E_0$, any $\{b^0_n\}_n \subset
\cup_\beta \alpha_1(\Cal H^0_\beta(\Cal N, \varphi))$ with
$\varphi(b^0_n)=0$, $\varphi({b_j^0}^*b^0_i)=\delta_{ij}$,
$\forall i,j,n$, and sp$(\{1\} \cup \{b^0_n\}_n)$ total in $\Cal
N$ gives an orthonormal basis of $\Cal B$ over $\Cal B_0$ by
letting $b_n=\alpha_1(b_n^0)$. To get such a set $\{b^0_n\}_n$
start with a total subset $1=a_0, a_1, a_2, ... \in \cup_\beta
\Cal H^0_\beta(\Cal N, \varphi)$ then apply the Gram-Schmidt
algorithm with respect to the scalar product $\langle \cdot, \cdot
\rangle_\varphi$.

We next show that $\{b_n\}_n$ this way defined automatically
satisfies $(3.3.1)$. By condition $(1.4.3), (a)$ and Kaplanski's
density theorem, there exist some finite, selfadjoint set of
unitary elements $v_k\in \Cal U_1= \{v\in \Cal U(\Cal
\alpha_1(\Cal N)) \mid v\Cal N v^* = \Cal N, {\text{\rm d}}\varphi
(v\cdot v^*) /\text{\rm d}\varphi, {\text{\rm d}}\varphi (v^*
\cdot v)/\text{\rm d}\varphi < \infty \}$ and scalars $c_k^j$ such
that $b'_j = \Sigma_k c_k^j v_k$ satisfy $\|b'_j\| \leq \|b_j\|$
and $\|b_j-b'_j\|_\varphi \leq \varepsilon(1+2{\text{\rm max}}_j
\|b_j\|)^{-1}$. Since $\|E_0(x_1x_2)\|_\varphi \leq
\|x_1x_2\|_\varphi \leq \|x_1\|\|x_2\|_\varphi$, $\forall x_1,x_2
\in \Cal B$, for $y \in \Cal N$ with $\|y\|\leq 1$ we get:
$$
\|E_0(b_i^*y\sigma_g(b_j))\|_\varphi \leq
\|E_0(b_i^*y\sigma_g(b'_j))\|_\varphi +
\|E_0(b_i^*y\sigma_g(b_j-b'_j))\|_\varphi
$$
$$
\leq \Sigma_k |c_k^j| \|E_0(b_i^*y\sigma_g(v_k))\|_\varphi +
\|b_i\| \|b_j-b^0_j\|_\varphi
$$
$$
\leq C\Sigma_k  \|E_0(b_i^*y\sigma_g(v_k))\|_\varphi +
\varepsilon/2,
$$
where $C = \text{\rm max} \{ |c^j_k| \}_{j,k}$. But if $v'_k =
\sigma_g(v_k)$ then $b_i^* y v_k' = b^*_i v'_k ({v'}^*_k y v'_k)$
and ${v'}_k^*yv'_k \in \Cal N$. Thus, $E_0(b_i^* y
v'_k)=\tilde{\varphi}(b_i^*v'_k){v'}^*_k y v'_k$. Thus, if we
denote $c={\text{\rm max}}\{ {\text{\rm d}}\varphi(v_k^* \cdot
v_k)/{\text{\rm d}}\varphi \}_k$ then $\|E_0(b_i^* y
v'_k)\|_\varphi \leq c |\tilde{\varphi}(b_i^*v'_k)|$ and from the
above estimates we get:

$$
\|E_0(b_i^*yb_j)\|_\varphi \leq cC\Sigma_k
|\tilde{\varphi}(b_i^*\sigma_g(v_k))| + \varepsilon/2.
$$

But $\tilde{\sigma}$ is mixing on $\tilde{\Cal N}$, in particular
on $\alpha_1(\Cal N)$. Thus, since $\tilde{\varphi}(b_i)=0$, we
get $\underset g \rightarrow \infty \to \lim
\tilde{\varphi}(b_i^*\sigma_g(v_k))=0$, $\forall i,k$, showing
that $(3.3.1)$ holds true.

Now, since for all $u \in \Cal U(P_0'\cap \tilde{M})$ we have
$ue_0u^* \in P_0'\cap \langle \Cal B, e_0 \rangle$ and $ue_0u^*$
is in the centralizer of $\phi$, we can apply Proposition 3.3 to
get $ue_0u^* \leq e_0$. By the faithfulness of $\phi$ this implies
$ue_0u^*=e_0$, or $ue_0=e_0u$. Applying this equality to the
separating vector $\hat{1}$ in $L^2(\tilde{\Cal M},
\tilde{\varphi})$, this gives $E_0(u)=u$. But since $E_0$ is
$\tilde{\varphi}$-preserving, it takes the centralizer of
$\tilde{\varphi}$ into the centralizer of $\varphi$ on $\Cal M$,
i.e., $E_0(\tilde{M})=M$. This yields $E_0(P_0'\cap
\tilde{M})\subset P_0'\cap M$, thus $u=E_0(u)\in P_0'\cap M$.
\hfill Q.E.D.

We end this Section with a technical lemma needed in the proof
of the main result in the next section.

\proclaim{3.5. Lemma} Let $Q \subset P$ be an inclusion of finite
von Neumann algebras and $q \in \Cal P(Q)$, $q'\in \Cal P(Q'\cap
P)$.

$1^\circ$. If $Q$ is quasi-regular in $P$ then $qQq'q$ is quasi-regular in
$qq'Pqq'$ (see $1.4.2$ in [Po3] for the definition
of quasi-regular subalgebras).

$2^\circ$. If $Q$ is regular in $P$ and $q\in \Cal P(Q)$ satisfies
$ctr_Q(q)=cz$, for some scalar $c$ and central projection $z\in
\Cal Z(Q)$, then $qQq$ is regular in $qPq$.
\endproclaim
\noindent
{\it Proof}. 1$^\circ$. We may clearly assume
$P$ has a normal faithful trace
$\tau$. It is then sufficient to prove that $\forall
x \in q\Cal N_P(Q)$ and $\varepsilon > 0$
$\exists z\in \Cal Z(Q)$
such that $\tau(1-z) \leq
\varepsilon$ and $qq'zxqq'z \in q\Cal N_{qq'Pqq'}(qQqq')$.
Let $x_1, \ldots, x_n \in P$ be so that
$Qx\in \Sigma_i x_i Q$ and $xQ \subset \Sigma_i Qx_i$.

For the given $\varepsilon > 0$ there exists $z \in \Cal Z(Q)$ and
finitely many partial isometries $v_1, v_2, \ldots, v_m \in Q$
such that $\tau(1-z) \leq \varepsilon$, $v_j^*v_j \leq q$ and
$\Sigma_j v_jv_j^* = z$. If we let $\{y_k\}_k$ be a relabelling of
the finite set $\{qq'zx_iv_jqq'z\}_{i,j} \cup
\{qq'zv_j^*x_iqq'z\}_{i,j}$, then we have
$$
(qq'Qqq')(qq'zxqq'z) \subset \Sigma_k y_k (qq'Qqq'),
$$
$$
(qq'zxqq'z)(qq'Qqq') \subset \Sigma_k (qq'Qqq')y_k.
$$

$2^\circ$. If $q=z \in \Cal Z(Q)$ then
$\forall u \in \Cal N_P(Q)$,
$v=zuz$ follows a partial isometry with $vv^*, v^*v \in \Cal Z(Q)$ and
$vQv^* = vv^*Qvv^*$. The proof of (2.1 in [JPo])
shows that $v$ can be extended to a unitary element in $zPz$
normalizing $Qz$. Thus, $z\Cal N_P(Q)z \subset \Cal N_{zPz}(Qz)''$,
implying that $zPz=\overline{\text{\rm sp}}^w z\Cal N_P(Q)z
\subset \Cal N_{zPz}(Qz)''$, i.e., $\Cal N_{zPz}(Qz)''=zPz$.

If $ctr_Q(q)=cz$ for some scalar $c$ and $z \in \Cal Z(Q)$,
then by the first part $Qz$ is regular in $zPz$. This reduces
the general case to the case $q \in P$ has scalar central trace
in $Q$. But then, if $u \in \Cal N_P(Q)$ we have $ctr(uqu^*)
=ctr(q)$ so there exists $v \in \Cal U(Q)$ such that $vuqu^*v^* = q$,
implying that $q(vu)q\in qPq$ is a unitary element in the normalizer
$\Cal N_{qPq}(qQq)$. Thus, $\text{\rm sp} Q \Cal N_{qPq}(qQq)Q
\supset {\text{\rm sp}} \Cal N_P(Q)$.

This
yields ${\text{\rm sp}} qQq \Cal N_{qPq}(qQq) qQq \supset q\Cal N_P(Q)q$
and since the right hand  term generates $qPq$ while the left hand
one is generated by $\Cal N_{qPq}(qQq)$, we get
$\Cal N_{qPq}(qQq)'' = qPq$.

\hfill Q.E.D.

\heading 4. Rigid embeddings into $N \rtimes_\sigma G$ are
absorbed by $L(G)$
\endheading

In this section we prove a key rigidity  result for inclusions of
algebras of the form $L(G) \subset M=N \rtimes_\sigma G$, in the
case $\sigma: G \rightarrow {\text{\rm Aut}}(N, \tau)$ is a
malleable mixing action. Thus, we show that if $Q \subset M$ is a
diffuse, relatively rigid von Neumann subalgebra whose normalizer
in $M$ generates a factor $P$ (see Section 4 in [Po3] for the
definition of relatively rigid subalgebras), then $Q$ and $P$ are
``absorbed'' by $L(G)$, via automorphisms of $M^\infty$ coming
from given gauged extensions of $\sigma$.

The proof follows an idea from ([Po1]): Due to malleability, the
algebra $M = N \rtimes_\sigma G$ can be perturbed continuously via
the gauge action, leaving only $L(G)$ fixed, thus forcing any
relatively rigid subalgebra $Q$ of $M$ to sit inside $L(G)$
(modulo some unitary conjugacy). The actual details of this
argument will require the technical results from the previous
sections. In particular, in order to apply Corollary 3.3 we'll
need $P$ not embeddable into $N$.

An example when $P$ is not embeddable into $N$ is when $Q$ is
already rigid in $P$ and $N$ has Haagerup's compact approximation
property (cf. 5.4.1$^\circ$ in [Po3]; for the definition of
Haagerup's  property for algebras see [Cho], or  2.0.2 in [Po3]).
In particular this is the case if $N$ is approximately finite
dimensional (AFD), or if $N=L(\Bbb F_n)$ for some $2 \leq n \leq
\infty$. Another example of this situation is when $N$ is abelian.

\proclaim{4.1. Theorem} Let $M$ be a factor of the form $M = N
\rtimes_\sigma G$, for some malleable mixing action $\sigma$ of a
discrete ICC group $G$ on a finite von Neumann algebra $(N,
\tau)$. Assume $Q\subset M^s$ is a diffuse, relatively rigid von
Neumann subalgebra such that $P=\Cal N_{M^s}(Q)''$ is a factor.
Assume also that no corner of $P$ can be embedded (non-unitally)
into $N$. If $\tilde{\sigma}$ is a gauged extension for $\sigma$,
then there exist a unique $\beta \in H(\tilde{\sigma})$ and a
unique $\theta_\beta\in {\text{\rm Aut}}_{\beta} (M;
\tilde{\sigma})$ such that the isomorphism $\theta_\beta: M^s
\simeq M^{s\beta}$ satisfies $\theta_\beta(P)\subset
L(G)^{s\beta}$.
\endproclaim
N.B.: The uniqueness of $\theta_\beta$ is modulo perturbations
from the left by inner automorphisms implemented by unitaries from
$L(G)^{s\beta}$.
\vskip .05in
\noindent
{\it Proof}. The
uniqueness is trivial by Theorem 3.1. We split the proof of the
existence into seven {\it Steps}. For the first six {\it Steps},
we assume $s=1$. Then in {\it Step} 7, we use the case $s=1$ to
settle the general case.

{\it Step 1}. $\exists \delta > 0$
such that $\forall t > 0, t\leq \delta$,
$\exists w(t) \in \tilde{M}, w(t)\neq 0$, satisfying  $w(t)y=\alpha_{t}(y)w(t),
\forall y\in Q$.

Let $\tilde{\sigma} : G \rightarrow
{\text{\rm Aut}} (\Cal N \subset \tilde{\Cal N}, \tilde{\varphi})$
be the given
gauged extension with gauge
$\alpha: \Bbb R \rightarrow {\text{\rm Aut}}
(\tilde{\Cal N}, \tilde{\varphi})$. With
the notations in 1.7,
$\alpha$ implements a continuous action
of $\Bbb R$ on the discrete decomposition
$(\tilde{\Cal M}, \tilde{\varphi})$.
In particular, this action implements a continuous action of
$\Bbb R$ on the type II$_1$ factor $\tilde{M}$, still denoted $\alpha$.
Thus, $\underset t \rightarrow 0 \to \lim \|\alpha_t(x)-x\|_2=0$,
$\forall x\in \tilde{M}$.

Since $Q \subset M$ is rigid, $Q \subset \tilde{M}$ is also rigid,
so that there exists $\delta > 0$ such that if $|t| \leq \delta$
then $\|\alpha_t(u)-u\|_2 \leq 1/2, \forall u \in \Cal U(Q).$ Let
$a(t)$ be the unique element of minimal norm-$\|\quad\|_2$ in
$\overline{\text{\rm co}}^w \{\alpha_t(u)u^*\mid u\in Q\}$. Since
$\|\alpha_t(u)u^*-1\|_2 \leq 1/2, \forall u \in \Cal U(Q)$, we
have $\|a(t)-1\|_2 \leq 1/2$, thus $a(t)\neq 0$.

By the uniqueness of
$a(t)$ we have
$\alpha_t(u)a(t)u^*=a(t), \forall u\in \Cal U(Q)$.
Thus, $a(t)y=\alpha_t(y)a(t), \forall y\in Q$,
which implies $a(t)^*a(t) \in Q'\cap \tilde{M}$ and
$a(t)a(t)^* \in \alpha_t(Q)'\cap \tilde{M}$.

Thus, if we denote by $w(t)$ the partial
isometry in the polar decomposition of $a(t)$,
$w(t)=a(t)(|a(t)|)^{-1}$,  then $w(t)\neq 0$ and
$w(t)y=\alpha_t(y)w(t)$.

{\it Step 2}. If $t, w=w(t)\neq 0$ are such that $wy=\alpha_t(y)w,
\forall y\in Q$, then there exists a partial isometry $w' \in
\tilde{M}$ such that $w'y=\alpha_{2t}(y)w', \forall y \in Q,$ and
$\tau(w'{w'}^*) > \tau(ww^*)^2/2$.

To prove this note first that if $v \in \tilde{M}$ is a unitary
element normalizing $Q$ and if $\sigma_v$ denotes the automorphism
$v\cdot v^*$ on $Q$, then $\alpha_t(v)wv^*$ satisfies
$$
(\alpha_t(v)wv^*)y = \alpha_t(v)w \sigma_{v^*}(y)v^*
$$
$$
=\alpha_t(v)\alpha_t(\sigma_{v^*}(y))w v^* = \alpha_t(v\sigma_{v^*}(y))wv^*
$$
$$
=\alpha_t(v\sigma_{v^*}(y)v^*)\alpha_t(v)wv^* = \alpha_t(y)(\alpha_t(v)wv^*).
$$

We claim there exists $v$ in the normalizer
$\Cal N_{\tilde{M}}(Q)$ of $Q$ in $\tilde{M}$ such that
$$
\tau(vww^*v^*\alpha_t^{-1}(w^*w)) > \tau(ww^*)^2/2.
$$
Indeed, for if we would have $\tau(vww^*v^*\alpha_t^{-1}(w^*w))
\leq  \tau(ww^*)^2/2$, $\forall v\in \Cal N_{\tilde{M}}(Q)$, then
the element $h$ of minimal $\|\quad\|_2$ in $\overline{\text{\rm
co}}^w \{v(ww^*)v^* \mid v\in \Cal N_{\tilde{M}}(Q)\}$ would
satisfy $0\leq h \leq 1$, $h \in \Cal N_{\tilde{M}}(Q)'\cap
\tilde{M}$, $\tau(h)=\tau(ww^*)$ and $\tau(h\alpha_t^{-1}(w^*w))
\leq \tau(ww^*)^2/2$. But $\Cal N_{\tilde{M}}(Q)'\cap \tilde{M}
\subset P'\cap \tilde{M}$, and by Corollary 3.3 the latter equals
$P'\cap M$. Since by hypothesis one has $P'\cap M = \Bbb C$, this
shows that $h\in \Cal N_{\tilde{M}}(Q)'\cap \tilde{M} = \Bbb C$.
Thus we get
$$
\tau(ww^*)^2 = \tau(h)\tau(\alpha_t^{-1}(w^*w))
\leq \tau(ww^*)^2/2,
$$
a contradiction.

For such $v \in \Cal N_{\tilde{M}}(Q)$, let
$a=\alpha_t(\alpha_t(v)wv^*)w$ and note that $ay=\alpha_{2t}(y)a,
\forall y\in Q$. Moreover, the left support projection of $a$ has
trace $\geq \tau(vww^*v^*\alpha_t^{-1}(w^*w)) > \tau(ww^*)^2/2$.
Thus, if we take $w'$ to be the partial isometry in the polar
decomposition of $a$, $w'=a|a|^{-1}$, then $w'y=\alpha_{2t}(y)w',
\forall y\in Q$. Also, since $\|a\| \leq 1$, we have
$\tau(w'{w'}^*) > \tau(ww^*)^2/2$.

{\it Step 3}. There exists a non-zero partial isometry
$w_1 \in \tilde{M}$
such that $w_1y=\alpha_1(y)w_1, \forall y\in Q$.

To prove this, let first $w$ be a non-zero partial isometry in
$\tilde{M}$, satisfying $wy=\alpha_{2^{-n}}(y)w, \forall y\in Q$,
for some large $n \geq 1$, as given by {\it Step 1}. Set $v_0=w$.
By {\it Step 2} and induction, there exist partial isometries $v_k
\in \tilde{M}, k=0, 1, 2, ...,$ such that $v_k y =
\alpha_{2^{-n+k}}(y)v_k$, $\forall y\in Q$, and $\tau(v_kv^*_k)>
\tau(v_{k-1}v_{k-1}^*)^2/2, \forall k\geq 1$. Taking $w_1=v_n$, it
follows that $w_1y=\alpha_1(y)w_1, \forall y\in Q$ and
$\tau(w_1w_1^*)> \tau(ww^*)^{2^n}/2^{{2^n}-1}\neq 0$.

{\it Step 4}. With the notations in 1.8, there exists a positive, non-zero
element $b \in Q'\cap \langle \Cal M, e_1 \rangle$
such that $Tr(b) < \infty$.

Indeed, if $w_1 \in \tilde{M}$ is as
given by {\it Step 3}, then $w_1^*\tilde{e}_1w_1$ is a non-zero
positive element in
${\text{\rm sp}} \tilde{\Cal M}_1 \tilde{e}_1 \tilde{\Cal M}_1
\subset \langle \tilde{\Cal M}_1, \tilde{e}_1 \rangle$
that commutes with $Q$ and satisfies
$0\neq \tilde{\phi}(w_1^*\tilde{e}_1w_1) < \infty$.
Define $b=\Cal F(w^*_1\tilde{e}_1w_1)
\in \langle \Cal M, \tilde{e}_1 \rangle \simeq
\langle \Cal M, e_1 \rangle$. Then $b\neq 0$  and
$0 \leq b \leq 1$. Since $\Cal F$ is $\tilde{\phi}$-preserving,
$\phi(b)=\tilde{\phi}(b) \leq 1$.

{\it Step 5}. There exist projections $q\in Q, q'\in Q'\cap P$ and
a partial isometry $v_0 \in \Cal G\Cal V_\beta(\Cal M, \varphi)$
for some $\beta \in H(\tilde{\sigma})$, such that $v_0v_0^* = qq'$
and $v_0^*Pv_0 \subset L(G)$.

By $1^\circ \implies 4^\circ$ in Theorem 2.1, there exist
non-zero projections $q\in Q$,
$p\in L(G)$, an isomorphism $\psi$ of $qQq$ into
$pL(G)p$ and a non-zero partial isometry
$v_0\in \Cal G\Cal V_{\beta}(\Cal M, \varphi)$,
for some $\beta \in H(\tilde{\sigma})$,
such that $v_0v_0^* \in (qQq)'\cap qMq$, $v_0^*v_0
\in \psi(qQq)'\cap pMp$ and $xv_0=v_0\psi(x), \forall x\in qQq$.

Since $\psi(qQq)$ is a diffuse von Neumann subalgebra in $pL(G)p$,
by Theorem 3.1 it follows that $\psi(qQq)'\cap pMp \subset
pL(G)p$, showing that
$$
v_0^*Qv_0 = v_0^* q Qqq' v_0=\psi(Q)v_0^*v_0 \subset L(G).
$$
But since $v_0v_0^* \in (qQq)'\cap qMq$, it follows that $v_0v_0^*
= qq'$ for some $q'\in Q'\cap M \subset \Cal N_M(Q)'' =P$. Thus,
$qq'\in P$ as well. By Lemma 3.7, $qQq'q$ is quasi-regular in
$qq'Pqq'$, implying that $v_0^*Qv_0$ is quasi-regular in
$v_0^*Pv_0$. By Theorem 3.1, this shows that $v_0^*Pv_0 \subset
L(G)$.

Note that, since $P'\cap M=\Bbb C$ and since we did not use up to
now the condition that $G$ is ICC, the rest of the conditions in
the hypothesis of 4.1 are sufficient to imply that $G$ has finite
radical (see Theorem 4.4 below).

{\it Step 6}. End of the proof of the case $s=1$.

With the notations in {\it Step 5}, since Ad$v_0^* \in {\text{\rm
Aut}}_\beta (M;\tilde{\sigma})$  and since $P$ and $L(G)$ are
factors, it follows that there exists an appropriate amplification
$\theta_\beta: M \simeq M^\beta$ of Ad$v_0^*$ such that
$\theta_\beta(P)\subset L(G)^\beta$.

{\it Step 7}. Proof of the general case.

For general $s$,
note first that $P=\Cal N_M(Q)''$
being a factor and $Q$ being diffuse, for any $1 \geq t > 0$
there exists a projection in $q\in Q$
of trace
$\tau(q)=t$ and either $q\in \Cal Z(Q)$ (in case
$Q$ is type I homogeneous) or $ctr_Q(q)=c1$ (in case $Q$ is of type II$_1$).

Thus, by replacing $P$ by a factor $P_0$ of the form $M_{n\times
n}(qPq)\simeq P^{1/s} \subset (M^s)^{1/s}=M$, for some $t, n$ with
$tn=s^{-1}$, and $Q$ by its subalgebra $Q_0=D_n \otimes qQq$,
where $D_n \subset M_{n\times n}(\Bbb C)$ is the diagonal
subalgebra, by Lemma 3.5 we get a subfactor $P_0$ of $M$ with $Q_0
\subset P_0$ a diffuse von Neumann subalgebra such that $\Cal
N_M(Q_0)''=P_0$, while $Q_0 \subset M$ still a rigid inclusion
(the latter due to 4.4 and 4.5 in [Po3]). The first part applies
to get $\theta_\beta \in {\text{\rm Aut}}_\beta(M;\tilde{\sigma})$
such that $\theta_\beta(P_0)\subset L(G)^\beta$ and since
$(P\subset M^s)=(P_0\subset M)^s$, an appropriate
$s$-amplification of $\theta_\beta$ carries $P$ into
$L(G)^{s\beta}$.

\hfill Q.E.D

Note that for commutative Bernoulli shifts (1.6.1) we could only
prove malleability in the case the base space $(Y_0, \nu_0)$ has
no atoms. To prove that 4.1 holds true for all commutative
Bernoulli shifts $\sigma$ and under much weaker conditions on $Q$,
we consider the following ``malleability-type'' condition:

\vskip .05in \noindent {\bf 4.2. Definition}.  Let $(N_1, \tau_1)$
be a diffuse abelian von Neumann algebra and $\sigma_1$ an action
of $G$ on $(N_1, \tau_1)$. $\sigma_1$ is {\it sub malleable}
(resp. {\it sub  s-malleable}) if it can be extended to a
malleable (resp. s-malleable) action $\sigma$ of $G$ on a larger
abelian von Neumann algebra $(N, \tau)$ such that there exists an
orthonormal basis $\{1\} \cup \{b_i\}_i \subset N$ of $N$ over
$N_1$ satisfying

$$
\underset g \rightarrow \infty \to \lim ({\text{\rm sup}}\{ \|
E^{N}_{N_1}(b^*_i y \sigma(g)(b_j))\|_\varphi \mid y \in N_1,
\|y\| \leq 1 \}) = 0, \forall i,j  \tag 4.2.1
$$
Recall that by $(1.6.1)$, classical Bernoulli $G$-actions with
non-atomic (diffuse) base space are s-malleable mixing. We next
show that Bernoulli $G$-actions with arbitrary base are sub
s-malleable mixing:

\proclaim{4.3. Lemma} Let $\sigma_1$ be the Bernoulli shift action
of $G$ on $(X, \mu) = \Pi_g (Y_0, \nu_0)_g$, where $(Y_0, \nu_0)$
is an arbitrary (possibly atomic) non-trivial standard probability
space. Then the action $\sigma_1$ it induces on $L^\infty(X, \mu)$
is sub s-malleable mixing. \endproclaim \vskip .05in \noindent
{\it Proof}.  Denote $A_0^0 = L^\infty(Y_0, \nu_0)$ and consider
the embedding $A_0^0 \subset A^0_0 \overline{\otimes}
L^\infty(\Bbb T, \lambda) \simeq L^\infty(\Bbb T, \lambda)=A^0$.
If $z \in L^\infty(\Bbb T, \lambda)$ is the Haar generating
unitary and $u = 1 \otimes z \in A^0$ then $\{u^n\}_n$ is an
orthonormal basis of $A^0$ over $A_0^0$. Let $N_1 =
\overline{\otimes}_g (A^0_0)_g$, $N = \overline{\otimes}_g
(A^0)_g$ and denote $\{b_n\}_n \subset N$ the set of elements with
$b_n = \otimes (u^{n_g})_g$, $n_g \in \Bbb Z$ all but finitely
many equal to $0$.

It is immediate to see that $\{b_n\}_n$ is an orthonormal basis of
$N$ over $N_1$ that checks condition $(4.2.1)$ with respect to the
Bernoulli shift action $\sigma$ of $G$ on $N=\overline{\otimes}_g
(A^0)_g$. Since $\sigma$ extends the Bernoulli shift action
$\sigma_1$ of $G$ on  $N_1 = \overline{\otimes}_g (A^0_0)_g$ and
since by $(1.6.1)$ it has mixing graded gauged extensions,
$\sigma_1$ follows sub s-malleable mixing. \hfill Q.E.D.

\proclaim{4.4. Theorem} Let $M_1$ be a factor of the form $M_1 =
N_1 \rtimes_{\sigma_1} G$, for some free action $\sigma_1$ of a
discrete group $G$ on an abelian von Neumann algebra $(N_1,
\tau_1)$. Let $Q\subset M_1^s$ be a diffuse, relatively rigid von
Neumann subalgebra and denote $P_1=\Cal N(Q)''$ the von Neumann
algebra generated by its normalizer in $M_1^s$.

$(i)$. If $\sigma_1$ is sub malleable and $P_1$ has atomic center
then $G$ has finite radical and for any minimal projection $p_1$
in the center of $P_1$ there exists a minimal projection $p$ in
the center of $L(G)$ and a unitary element $u\in M_1^s$ such that
$u(P_1p_1)^{t_1} u^* \subset (L(G)p)^{st}$, where
$t_1=\tau(p_1)^{-1}$, $t=\tau(p)^{-1}$. If moreover $P_1$ is a
factor and $G$ is ICC then there exists $u \in \Cal U(M_1^s)$ such
that $uP_1u^* \subset L(G)^s$. Also, $u$ is unique with the above
property, modulo perturbation from the left by a unitary in
$L(G)^s$.

$(ii)$. If $\sigma_1$ is sub s-malleable, $Q$ is of type
${\text{\rm II}}_1$ and $G$ is ICC then there exists $u \in \Cal
U(M_1^s)$ such that $uP_1u^* \subset L(G)^s$, unique modulo
perturbation from the left by a unitary in $L(G)^s$.
\endproclaim
\vskip .05in \noindent {\it Proof}. Let $\sigma: G \rightarrow
{\text{\rm Aut}}(N,\tau)$ be a malleable mixing extension of
$\sigma_1$ which satisfies $(4.2.1)$ and which is either
malleable, in the case $(i)$, or s-malleable, in the case $(ii)$.
Denote $M = N \rtimes_\sigma G$. We first show that in both cases
we have $P'\cap M^s = P'\cap M_1^s$. One has in fact the following
more general result:

\proclaim{4.5. Lemma} Let $(N_1, \tau_1) \subset (N, \tau)$ be an
embedding of finite von Neumann algebras for which there exists an
orthonormal basis $\{1\}\cup \{b_n\}_n$ of $N$ over $N_1$
satisfying condition $(4.2.1)$. Let $\sigma: G \rightarrow
{\text{\rm Aut}}(N, \tau)$ be a properly outer mixing action
leaving $N_1$ globally invariant. Denote $M_1 = N_1 \rtimes G$,
regarded as a von Neumann subalgebra of $M=N\rtimes G$. Assume
$P_1 \subset M_1$ is so that no corner of $P_1$ can be embedded
into $N$. If $x\in M$ satisfies $P_1 x \subset \Sigma_i x_i M_1$,
for some finite set $x_1, x_2, \ldots, x_n \in M$, then $x\in
M_1$. In particular, $P_1'\cap M= P_1'\cap M_1$.
\endproclaim
\vskip .05in \noindent {\it Proof}. If we put $\Cal T_0 = N_1,
\Cal T = N$, $\Cal B_0 = M_1, \Cal B=M$ then condition $(4.2.1)$
on $\{b_n\}_n$ shows that the hypothesis of Proposition 3.3 is
satisfied. Thus, if $E_0$ denotes the trace preserving conditional
expectation of $M$ onto $M_1$ and $M_1 \overset{E_0}\to\subset M
\subset \langle M,e_0 \rangle$ the corresponding basic
construction then any $a \in P_1'\cap \langle M, e_0 \rangle$ with
$a \geq 0, Tr(a) < \infty$ satisfies $e_0ae_0=a$.

But if $x\in M$ satisfies the condition in the hypothesis then the
orthogonal projection $f_x$ of $L^2(M)$ onto the closure $\Cal H$
of sp$P_1xM_1$ in $L^2(M)$ then $f_x \in P_1'\cap \langle M, e_0
\rangle$ (because $P\Cal HM_1 \subset \Cal H$) and dim$\Cal
H_{M_1} \leq n <\infty$ (because $\Cal H$ is contained in the
right $M_1$-Hilbert module $\overline{\Sigma_j x_j M_1}\subset
L^2(M)$). Thus, $Tr(f_x) < \infty$, implying that $f_x \leq e_0$,
equivalently $\Cal H \subset L^2(M_1)$. In particular $x \in
L^2(M_1)$, thus $x \in M_1$. \hfill Q.E.D.

\vskip .05in \noindent {\it Proof of} $(i)$.  $N$ being abelian
and $P_1$ of type II$_1$, it follows that no corner of $P_1$ can
be embedded into $N$. By 4.5 this implies $P_1'\cap M^s = P_1'\cap
M_1^s=\Bbb C$.

By $3.5.2^\circ$ a suitable amplification $P$ of the factor
$P_1p_1$ is a unitally embedded into $M$ and contains a diffuse
subalgebra $Q\subset P$ such that $Q \subset M$ is rigid and $\Cal
N_M(Q)''=P$. Since $\sigma$ is malleable mixing and has gauged
extension with trivial spectrum (the gauged extension is even
abelian), by the first 5 Steps of the proof of 4.2 and the remark
at the end of Step 5, $G$ has finite radical and there exists a
non-zero partial isometry $v_0 \in M$ such that $v_0^*v_0 \in P$,
$v_0Pv_0^* \subset L(G)\subset M_1$. A suitable amplification
$u\in M^s$ of $v_0$, will then satisfy the conjugacy condition and
by 4.5 $u \in M_1^s$. The uniqueness is clear by 4.5.

\vskip .05in \noindent {\it Proof of} $(ii)$. Let $\tilde{\sigma}:
G \rightarrow \text{\rm Aut}(\tilde{N}, \tilde{\tau})$, $\alpha:
\Bbb R \rightarrow \text{\rm Aut}(\tilde{N}, \tilde{\tau})$,
$\beta \in \text{\rm Aut}(\tilde{N}, \tilde{\tau})$, $\beta^2=id$,
give a graded gauged extension for $\sigma$. Denote
$\tilde{M}=\tilde{N} \rtimes_{\tilde{\sigma}} G$.

We first prove that there exists a non-zero partial isometry $w\in
\tilde{M}$ such that $w^*w\in Q'\cap M$, $ww^* \in \alpha_1(Q'\cap
M)$, $wy=\alpha_1(y)w, \forall y\in Q$.

To this end, note first that Step 1 of the proof of Theorem 4.1,
which only used the fact that $\sigma$ is malleable, shows that
for all $t$ sufficiently small there exists a non-zero partial
isometry $v=v(t)\in \tilde{M}$ such that $v^*v\in Q'\cap
\tilde{M}$, $vv^* \in \alpha_t(Q)'\cap \tilde{M}$,
$vy=\alpha_t(y)v, \forall y \in Q$. But since $Q$ is of type
II$_1$, no corner of $Q$ can be embedded into $N$. Thus, by
Theorem 3.2 we have $Q'\cap \tilde{M}=Q'\cap M$ and so we get:
$$
v^*v\in Q'\cap M, vv^* \in \alpha_t(Q'\cap M), vy=\alpha_t(y)v,
\forall y \in Q. \tag 4.4.1
$$

Assuming now that for some $0 < t <1$ there exists a partial
isometry $v\in \tilde{M}$ such that $v^*v \in Q'\cap M$, $vv^* \in
\alpha_t(Q'\cap M)$ and $vy = \alpha_t(y)v, \forall y\in Q$, we
show that there exists a partial isometry $v'\in \tilde{M}$
satisfying $\|v'\|_2=\|v\|_2$, ${v'}^*v' \in Q'\cap M$, $v'{v'}^*
\in \alpha_t(Q'\cap M)$ and $v'y = \alpha_{2t}(y)v', \forall y\in
Q$. This will of course prove the existence of $w$, by starting
with some $t=2^{-n}$ for $n$ sufficiently large, then proceeding
by induction until we reach $t=1$.

Applying $\beta$ to $vy=\alpha_t(y)v$ and using that $\beta(x)=x,
\forall x\in M$ and $\beta \alpha_t = \alpha_{-t} \beta$, we get
$\beta(v)y=\alpha_{-t}(y)\beta(v)$, $\forall y\in Q$. Taking
adjoints and plugging in $y^*$ for $y$ we further get $y\beta(v^*)
= \beta(v^*) \alpha_{-t}(y), \forall y\in Q$. Thus:
$$
\alpha_t(y)v\beta(v^*) = v y \beta(v^*) = v \beta(v^*)
\alpha_{-t}(y), \forall y\in Q. \tag 4.4.2
$$
If we now apply $\alpha_t$ to the last and first term of $(4.4.2)$
and denote $v'=\alpha_t(v \beta(v^*))$, then we get:
$$
v'y = \alpha_{2t}(y)v', \forall y\in Q. \tag 4.4.3
$$
Moreover, since $v^*v\in Q \subset M$, we have $\beta(v^*v)=v^*v$
so $v\beta(v^*)$ is a partial isometry and it has the same range
as $v$. Thus $v'$ is a partial isometry and $\|v'\|_2 = \|v\|_2$.
Since $v'$ is an intertwiner, ${v'}^*v'\in Q'\cap \tilde{M}=Q'\cap
M$ and similarly $v'{v'}^* \in \alpha_{2t}(Q'\cap M)$.

With the non-zero partial isometry $w \in \tilde{M}$ satisfying
$wy=\alpha_1(y)w, \forall y\in Q$ we can apply  Steps 4 and 5 in
the proof of 4.2 to get a non-zero partial isometry $v_0 \in M$
such that $v_0^*v_0 \in Q'\cap M$ and $v_0Qv_0^* \subset L(G)$.

To end the proof we use a maximality argument. Thus, we consider
the set $\Cal W$ of all families $(\{p_i\}_i, u)$ where
$\{p_i\}_i$ are partitions of 1 with projections in $Q'\cap M$,
$u\in M$ is a partial isometry with $u^*u=\Sigma_i p_i$ and
$u(\Sigma_i Qp_i)u^* \subset L(G)$. We endow $\Cal W$ with the
order given by $(\{p_i\}_i, u) \leq (\{p'_j\}_j, u')$ if
$\{p_i\}_i \subset \{p'_j\}_j$, $u=u'(\Sigma_i p_i)$. $(\Cal W,
\leq)$ is clearly inductively ordered.

Let $(\{p_i\}_i, u)$ be a maximal element. If $u$ is a unitary
element, then we are done. If not, then denote $q'= 1 -\Sigma_i
p_i \in Q'\cap M$ and take $q\in Q$ such that $\tau(qq')=1/n$ for
some integer $n \geq 1$. Denote $Q_0=M_{n \times n}(qQqq')$
regarded as a von Neumann subalgebra of $M$, with the same unit as
$M$. By ([Po3]), it follows that $Q_0 \subset M$ is rigid. Thus,
by the first part there exists a non-zero partial isometry $w\in
M$ such that $w^*w\in Q_0'\cap M$ and $wQ_0w^* \subset L(G)$.
Since $qq'\in Q_0$ has scalar central trace in $Q_0$, it follows
that there exists a non-zero projection in $w^*w Q_0 w^*w$
majorised by $qq'$ in $Q_0$.

It follows that there exists a non-zero projection $q_0 \in
qq'Q_0qq' = qQqq'$ and a partial isometry $w_0 \in M$ such that
$w_0^*w_0 = q_0$ and $w_0(qQqq')w_0^* \subset L(G)$. Moreover, by
using the fact that $Q$ is diffuse, we may shrink $q_0$ if
necessary so that to be of the form $q_0 = q_1 q'\neq 0$ with $q_1
\in \Cal P(Q)$ of central trace equal to $m^{-1} z$ for some $z\in
\Cal Z(Q)$ and $m$ an integer. But then $w_0$ trivially extends to
a partial isometry $w_1 \in M$ with $w_1^*w_1 = q'z\in Q'\cap M$
and $w_1 Qw_1^* \subset L(G)$. Moreover, since $L(G)$ is a factor,
we can multiply $w_1$  from the left with a unitary element in
$L(G)$ so that $w_1w_1^*$ is perpendicular to $uu^*$. But then
$(\{p_i\}_i \cup \{q'z\}, u_1)$, where $u_1=u+w_1$, is clearly in
$\Cal W$ and is (strictly) larger than the maximal element
$(\{p_i\}_i, u)$, a contradiction. \hfill Q.E.D.

\vskip .05in Theorems 4.1, 4.4 allow in fact the control of
inclusions $\theta_\beta(P) \subset L(G)^{s\beta}$ from above as
well:

\proclaim{4.6. Corollary} Let $G$ be a discrete ICC group, $\sigma
: G \rightarrow {\text{\rm Aut}}(N, \tau)$, $M=N\rtimes_\sigma G$,
$s > 0$ and $Q\subset M^s$ a diffuse, relatively rigid von Neumann
subalgebra, with $P=\Cal N_{M^s}(Q)''$ a factor. Assume that
either $\sigma$ is malleable mixing with $P$ not embeddable into
$N$, or that $N$ is abelian and $\sigma$ is sub malleable mixing.
If $P \supset L(G)^s$ then $P = L(G)^s$, and thus $Q \subset
L(G)^s$ as well.
\endproclaim
\vskip .05in
\noindent
{\it Proof}. Since $P, L(G)^s$ are factors, it is sufficient to prove that
there exists a non-zero projection $p\in \Cal P(L(G)^s) \subset \Cal P(P)$
such that $pPp \subset L(G)$.

In case  $\sigma$ is malleable mixing, by Theorem 4.1 it follows
that if $\tilde{\sigma}$ is a gauged extension for $\sigma$ then
there exists $\beta \in S(\tilde{\sigma})$ and $\theta_\beta \in
{\text{\rm Aut}}_\beta (M;\tilde{\sigma})$ such that $\theta_\beta
(P)\subset L(G)^{s\beta}$. In particular, since $P\supset L(G)^s$
this implies $\theta_\beta(L(G)^s) \subset L(G)^{s\beta}$. But
then Theorem 3.1 implies $\beta = 1$ and that $\theta_\beta$ is
actually implemented by a unitary element in $L(G)^s$. The case
$N$ abelian, $\sigma$ sub-malleable follows similarly, using 4.4
(a) instead of 4.1. \hfill Q.E.D.

\heading 5. Strong rigidity of the inclusions $L(G)\subset N
\rtimes_\sigma G$
\endheading

The ``absorbtion'' results in the previous Section allow us to
prove here that any isomorphism between amplifications of factors
of the form $N \rtimes_\sigma G$, with $G$ satisfying a weak
rigidity property, $N$ approximately finite dimensional (or merely
having Haagerup's property) and $\sigma$ malleable mixing, can be
perturbed by an automorphism coming from a given gauged extension
so that to carry the subalgebras $L(G)$ onto each other. We
consider the following terminology:

\vskip .05in \noindent {\it 5.1. Definitions}. $1^\circ$. Given a
group $\Gamma$, an infinite subgroup $\Lambda\subset \Gamma$ is
{\it wq-normal} in $\Gamma$ if for any intermediate subgroup
$\Lambda \subset H \varsubsetneq  \Gamma$ there exists $g \in
\Gamma \setminus H$ such that  $gHg^{-1} \cap H$ is infinite. In
particular, if there exist subgroups $\Lambda=\Lambda_0 \subset
\Lambda_1 \subset ... \subset \Lambda_n=\Gamma$ with $\Lambda_i
\subset \Lambda_{i+1}$ normal $\forall 0\leq i \leq n-1$, then
$\Lambda$ is wq-normal in $\Gamma$. An example of a non wq-normal
inclusion of groups is $\Lambda \subset \Gamma=\Lambda
* K_0$ with $K_0$ non trivial. However, $\Lambda \subset (\Lambda
* K_0) \times K$ is a wq-normal inclusion
whenever $\Lambda, K$ are infinite groups.

$2^\circ$. A discrete group $\Gamma$ is {\it w-rigid} if it
contains an infinite normal subgroup with the relative property
(T) of Kazhdan-Margulis ([Ma]; see also [dHV])). Also, we denote
by $w\Cal T_0$ the class of groups $\Gamma$ having a non virtually
abelian wq-normal subgroup $\Lambda \subset \Gamma$ with the
relative property (T).

\vskip .05in If an infinite group $H$ has the property (T) of
Kazhdan ([Ka]) then any group $G$ having $H$ as a normal subgroup
is w-rigid. For instance $G=H\rtimes K$ for some arbitrary group
$K$ acting on $H$ by (possibly trivial) automorphisms. Other
examples of w-rigid groups are the groups $G=\Bbb Z^2 \rtimes
\Gamma$ with $\Gamma \subset SL(2, \Bbb Z)$ non-amenable (cf.
[Ka], [Ma], [Bu]), and $G=\Bbb Z^m \rtimes \Gamma$ with $\Gamma$
arithmetic lattice in $SO(n,1)$ or $SU(n,1)$, for suitable $m$ and
suitable actions of $\Gamma$ on $\Bbb Z^m$ (cf. [Va]). Also, note
that the class $w\Cal T_0$ is closed to normal and finite index
extensions and to inductive limits.

\proclaim{5.2. Theorem} Let $G_i$ be w-rigid ICC groups,
$\sigma_i: G_i \rightarrow {\text{\rm Aut}}(N_i, \tau_i)$
malleable mixing actions with gauged extensions $\tilde{\sigma}_i$
and $M_i=N_i \rtimes_{\sigma_i} G_i$, $i=0,1$. Let $\theta : M_0
\simeq M_1^{s}$ be an isomorphism, for some $s > 0$. Then there
exist unique $\beta_i \in S(\tilde{\sigma}_i)$ and unique
$\theta^i_{\beta_i} \in {\text{\rm Aut}}_{\beta_i}(M_i;
\tilde{\sigma}_i)$ such that $\theta^1_{\beta_1}(\theta
(L(G_0)))=L(G_1)^{s\beta_1}$,
$\theta(\theta^0_{\beta_0}(L(G_0)))=L(G_1)^{s\beta_0}$. Moreover,
$\beta_0=\beta_1$.
\endproclaim
N.B. The uniqueness above is modulo perturbations by inner
automorphisms implemented by unitaries from the appropriate
amplifications of $L(G_i)$. \vskip .05in \noindent {\it Proof}.
Let $H_i \subset G_i$ be infinite wq-normal, relatively rigid
subgroups, $i=0,1$. Thus, $Q_i=L(H_i)$ is diffuse and the von
Neumann algebra $P_i$ generated by its normalizer $\Cal
N_{M_i}(Q_i)$ in $M_i$ contains $L(G_i)$ (because $H_i$ is normal
in $G_i$). Moreover, by Theorem 3.1, $P_i$ is contained in
$L(G_i)$. Thus, $P_i=L(G_i)$.

By applying Theorem 4.1 to $\theta(Q_0)\subset \theta(P_0) \subset
M_1^s$, it follows that there exist unique $\beta_1 \in
S(\tilde{\sigma}_1)$ and $\theta^1_{\beta_1} \in {\text{\rm
Aut}}_{\beta_1}(M_1;\tilde{\sigma}_1)$ such that
$\theta_{\beta_1}^1( \theta (P_0))\subset L(G_1)^{s\beta_1}$.
Thus, $(\theta^1_{\beta_1}\circ \theta)^{-1}(L(G_1)) \supset P_0 =
L(G_0)$ and by Corollary 4.6 we actually have equality. \hfill
Q.E.D.

We mention in separate statements the case of actions on
probability spaces, where due to Theorem 4.4 we can require the
actions $\sigma$ to be merely sub malleable (resp. sub
s-malleable). The proof is identical to the one for 5.2 above,
using $4.4.(i)$ (resp. $(4.4.(ii)$) in lieu of 4.1, and it is thus
omitted.

\proclaim{5.3. Theorem} Let $\sigma: G \rightarrow {\text{\rm
Aut}}(N, \tau)$, $\sigma_0: G_0 \rightarrow \text{\rm Aut}(N_0,
\tau_0)$ be free ergodic actions of discrete groups on abelian von
Neumann algebras $N,N_0$. Denote $M = N \rtimes_{\sigma} G$, $M_0
= N_0 \rtimes_{\sigma_0} G_0$. Let also $H \subset G_0'\subset
G_0$ be subgroups with $H$ wq-normal in $G_0'$. Assume:

$(a)$. $G$ is ICC; $\sigma$ is sub malleable mixing.

$(b)$. $H$ is w-rigid, $\{hgh^{-1} \mid h\in H\}$ is infinite
$\forall g \in G_0, g\neq e$; $\sigma_{0|H}$ is ergodic.

\vskip .05in

If $\theta:M_0 \simeq M^s$ is an isomorphism, then there exists
$u\in \Cal U(M^s)$ such that $u\theta(L(G_0'))u^* \subset L(G)^s$.
In particular, if $(\sigma, G)$ satisfies $(a)$ and $G_0$ is
w-rigid ICC then $u\theta(L(G_0))u^* \subset L(G)^s$. Moreover, if
both $G_0, G$ are w-rigid ICC and both $\sigma_0, \sigma$ sub
malleable mixing then $u\theta(L(G_0))u^* = L(G)^s$.
\endproclaim

\proclaim{5.3'. Theorem} Instead of $(a), (b)$ in $5.3$, assume

$(a')$. $G$ is ICC; $\sigma$ is sub s-malleable mixing.

$(b')$. $H$ has the relative property $(\text{\rm T})$ in $G_0$
$($i.e. $(G_0, H)$ is a property $(\text{\rm T})$ pair$)$ and is
not virtually abelian.

\vskip .05in

If $\theta:M_0 \simeq M^s$ is an isomorphism, then there exists
$u\in \Cal U(M^s)$ such that $u\theta(L(G_0'))u^* \subset L(G)^s$.
Moreover, if both $G_0, G$ are ICC and in the class $w\Cal T_0$
and both $\sigma_0, \sigma$ sub s-malleable mixing then
$u\theta(L(G_0))u^* = L(G)^s$.
\endproclaim

Theorems $4.1, 5.2, 5.3$, $5.3'$ are key ingredients in the proof
of strong rigidity results for isomorphisms of cross product
factors $N_0 \rtimes_{\sigma_0} G_0$, $N \rtimes_\sigma G$ for
$G_0$ w-rigid or in the class $w\Cal T_0$ and $\sigma$ commutative
or non-commutative Bernoulli shifts, in ([Po6]). This will allow
us to classify large classes of factors $N \rtimes_\sigma G$, with
explicit calculations of various invariants, such as $\mycal F(N
\rtimes_\sigma G)$.

In this paper we only mention a straightforward application of
Theorems 5.2 and of results from ([Po3]): It shows that the
fundamental group $\mycal F(N \rtimes_\sigma G)$ is equal to
$S(\tilde{\sigma})$ for any gauged extension $\tilde{\sigma}$ of a
malleable mixing $\sigma$, in the case $G=\Bbb Z^2 \rtimes
\Gamma$, where $\Gamma$ is a subgroup of finite index in $SL(2,
\Bbb Z)$, or if $G=\Bbb Z^N \rtimes \Gamma$ where $\Gamma$ is an
arithmetic lattice in either $SU(n, 1)$ or $SO(2n, 1)$ suitably
acting on $\Bbb Z^N$, or more generally, if $G$ is a finite
product of any of the above groups. Taking $\sigma$ to be
Connes-St\o rmer Bernoulli or Bogoliubov shifts, this result
already gives many examples of factors and equivalence relations
with arbitrarily prescribed countable fundamental group.

\proclaim{5.4. Corollary} $1^\circ$. Let $G$ be a w-rigid ICC
group and $\sigma$ a properly outer mixing action of $G$ on a AFD
von Neumann algebra $(N, \tau)$ $($more generally a finite von
Neumann algebra with Haagerup's approximation property$)$. Assume
$\sigma$ is malleable mixing $($resp. $N$ abelian and $\sigma$ sub
malleable mixing$)$. If $\mycal F(L(G))=\{1\}$ then $\mycal
F(N\rtimes_\sigma G)=S(\tilde{\sigma})$ for any gauged extension
$\tilde{\sigma}$ of $\sigma$ $($resp. $\mycal F(N\rtimes_\sigma
G)=\{1\})$. In particular, this is the case if $L(G)$ is a HT
factor and $\beta^{^{HT}}_n(L(G))\neq 0, \infty$ for some $n$.

$2^\circ$. Let $G_i$ be w-rigid ICC groups and $\sigma_i$
malleable mixing actions of $G_i$ on AFD algebras $(N_i, \tau_i)$,
$i=0,1$ $($resp. sub malleable mixing, with $N_i$ abelian,
$i=0,1)$. Assume $L(G_i)$ are HT factors with
$\beta^{^{HT}}_n(L(G_0))=0$ and $\beta^{^{HT}}_n(L(G_1))\neq
0,\infty$ for some $n$. Then $N_0\rtimes_{\sigma_0} G_0$ is not
stably isomorphic to $N_1\rtimes_{\sigma_1} G_1$.
\endproclaim
\noindent {\it Proof}. Assume first that $\sigma$ in part
$1^\circ$ and $\sigma_i$ in part 2$^\circ$ are malleable.

Under the hypothesis in 2$^\circ$, denote $M_i = N_i
\rtimes_{\sigma_i} G_i$ and assume $\theta : M_0 \simeq M_1^s$ is
an isomorphism, for some $s
> 0$. Let $\tilde{\sigma}_1$ be a gauged extension for $\sigma_1$.
By 5.2 there exists $\beta \in S(\tilde{\sigma}_1)$ and an
isomorphism $\theta_\beta : M_1^s \simeq M_1^{s\beta}$ such that
$\theta_\beta (\theta (L(G_0))) = L(G_1)^{s\beta}$.

Thus, if we take $G_i=G$, $N_i=N$, $\sigma_i=\sigma$ as in part
$1^\circ$, then $\mycal F(L(G))=\{1\}$ implies $st=1$, thus $s =
1/t \in S(\tilde{\sigma}_1)$, showing that $\mycal F(N
\rtimes_\sigma G) \subset S(\tilde{\sigma})$. Since we always have
$S(\tilde{\sigma})\subset \mycal F(N \rtimes_\sigma G)$, the
equality follows.

The rest of 1$^\circ$ and part 2$^\circ$ are now trivial, by the
first part of the proof and ([Po3]).

If we assume $N, N_i$ abelian and $\sigma, \sigma_i$, $i=0,1,$
sub-malleable, then the proofs are the same, but using 5.3 instead
of 5.2. \hfill Q.E.D.

\proclaim{5.5. Corollary} Let $S$ be a multiplicative subgroup in
$\Bbb R_+^*$. For each $ n \geq 2$ and $k \geq 1$ there exists a
properly outer action $\sigma_{n,k}$ of $(\Bbb F_n)^{k}=\Bbb F_n
\times \ldots \times \Bbb F_n$ $(k$ times$)$ on an AFD $\text{\rm
II}_1$ factor $N$ such that $\mycal F(N \rtimes_{\sigma_{n,k}}
\Bbb F^k_n) = S, \forall n,k,$ and such that for each $n \geq 2$
the factors $\{N \rtimes_{\sigma_{n,k}} (\Bbb F_n)^k\}_k$ are
non-stably isomorphic. Moreover, $N$ can be taken generated by at
most $|S|$ elements as a von Neumann algebra. In particular, if
$S$ is countable then $N$ can be taken the hyperfinite ${\text{\rm
II}}_1$ factor $R$.
\endproclaim
\noindent {\it Proof}. Let $G_{n,k}=(\Bbb Z^2)^k \rtimes (\Bbb
F_n)^k=(\Bbb Z^2 \rtimes \Bbb F_n)^k$, where the action of $\Bbb
F_n$ on $\Bbb Z^2$ comes from an embedding of finite index $\Bbb
F_n \subset SL(2, \Bbb Z)$. Let $\sigma_{n,k}'$ be a Connes-St\o
rmer Bernoulli shifts action of $G_{n,k}$ on an AFD II$_1$ factor
$N$, with gauged extension $\tilde{\sigma}_{n,k}'$ satisfying
$S(\tilde{\sigma}_{n,k}')=S$ (cf. last paragraph in 1.6.2). Since
$G_{n,k}$ are w-rigid, by 5.3 the factors $M_{n,k} = N
\rtimes_{\sigma_{n,k}'} G_{n,k}$ satisfy $\mycal F(M_{n,k})=S$.
Also, for each $n \geq 2$ $\{M_{n,k}\}_k$ are mutually non stably
isomorphic, because by [Po3] we have $\beta^{^{HT}}_k(L(G_{n,k}))=
\beta_k((\Bbb F_n)^k)\neq 0,\infty$ while
$\beta^{^{HT}}_j(L(G_{n,k}))=0$ for all $j \neq k$. But $N
\rtimes_{\sigma_{n,k}'} (\Bbb Z^2)^k \simeq R$, so denoting by
$\sigma_{n,k}$ the corresponding action of $(\Bbb F_n)^k$ on $N
\rtimes_{\sigma_{n,k}'} (\Bbb Z^2)^k$, the first part of the
statement follows. The last part is now trivial, since for $S$
countable we can take the base of the shift to be a separable type
I factor. \hfill Q.E.D.

\vskip .05in Note that in the case of Connes-St\o rmer
$(M_{k\times k}, \varphi_0)$-Bernoulli shift actions $\sigma : G
\rightarrow {\text{\rm Aut}}(N, \tau)$ (1.6.2), the cross product
factors $M=N \rtimes_\sigma G$ can be realized as a generalized
group measure space construction as in ([FM]), i.e., there exists
a standard ( = countable, ergodic, measure preserving) equivalence
relation $\Cal R$ on the standard probability space such that
$M=M(\Cal R)$. Indeed, this is because $\sigma$ normalizes the
``main diagonal'' Cartan subalgebra $A$ of $N$, thus normalizing
the corresponding hyperfinite equivalence relation $\Cal
R_{A\subset N}$ as well, inducing altogether an equivalence
relation $\Cal R=\Cal R_{A \subset M}$, with trivial cocycle.

Moreover, in each of these cases the automorphisms in
Aut$_\beta(M; \tilde{\sigma})$, given by the gauged extensions
$\tilde{\sigma}$ that come with the construction of $\sigma$, can
be chosen to normalize the Cartan subalgebra $A \subset N \subset
M$ as well. Thus, in each of these cases one has
$S(\tilde{\sigma}) \subset \mycal F(\Cal R)$. For non-commutative
Bernoulli shifts $\sigma$ this observation has been exploited in
[GeGo] to prove the existence of standard equivalence relations
with the fundamental group countable and containing a prescribed
countable set. Due to the obvious inclusion $\mycal F(\Cal R)
\subset \mycal F(M(\Cal R))$, Corollaries 5.3 and 5.4 allow us to
actually obtain precise calculations of fundamental groups of the
corresponding equivalence relations:

\proclaim{5.6. Corollary} Let $S$ be a countable multiplicative
subgroup in $\Bbb R_+^*$. There exist standard equivalence
relations $\Cal R$ such that $\mycal F(\Cal R)=S$. Moreover, given
any $ n \geq 2, k \geq 1$, $\Cal R$ can be taken of the form $\Cal
R = \Cal R_n \rtimes_{\sigma_{n,k}} (\Bbb F_n)^k$, where $\Cal R_n
\subset \Cal R$ is an ergodic  hyperfinite sub-equivalence
relation and $\sigma_{n,k}$ is a free ergodic  action of $(\Bbb
F_n)^k$ on the probability space, which normalizes $\Cal R_n$ and
acts outerly on it. Also, the equivalence relations $\{\Cal R_n
\rtimes (\Bbb F_n)^k\}_k$ can be taken mutually non-stably
isomorphic.
\endproclaim
\noindent {\it Proof}. In each of the examples used in the proof
of 5.4, for fixed $n \geq 2, k \geq 1,$ take $A$ to be the ``main
diagonal'' Cartan subalgebra of $N\simeq R$ which is normalized by
the action $\sigma_{n,k}'$ of $(\Bbb Z^2)^k \rtimes (\Bbb F_n)^k$.
Thus, the hyperfinite equivalence relation $\Cal R_{A\subset R}$
is invariant to $\sigma_{n,k}'$. It follows that if we let $\Cal
R_n$ be the equivalence relation implemented by $\Cal R_{A\subset
R}$ and ${\sigma'_{n,k}}_{|(\Bbb Z^2)^k}$, and $\Cal R$ be the
equivalence relation implemented jointly by $\Cal R_{A\subset R}$
and the action $\sigma'_{n,k}$ of $(\Bbb Z^2)^k \rtimes (\Bbb
F_n)^k$ then $M=R \rtimes_{\sigma'_{n,k}} ((\Bbb Z^2)^k \rtimes
(\Bbb F_n^k))$ coincides with the factor $M(\Cal R)$ associated
with $\Cal R$. Thus, $\mycal F(\Cal R) \subset \mycal F(M) = S$.
Since we also have $S \subset \mycal F(\Cal R)$, we are done.
\hfill Q.E.D.

\vskip .05in \noindent {\bf 5.7. Remarks}. $1^\circ$. The examples
of factors and equivalence relations with prescribed countable
fundamental group $S \subset \Bbb R_+^*$ can alternatively be
described as follows: Let $(Y_0,\nu_0)$ be an atomic probability
space with the property that the ratio set $\{\nu_0(x)/\nu_0(y)
\mid x,y\in Y_0 \}$ generates the group $S$. Let $G=\Bbb Z^2
\rtimes SL(2,\Bbb Z)$. Let $(X,\mu)=\Pi_g (Y_0,\nu_0)_g$. Denote
by $\Cal R_0$ the equivalence relation on $(X,\mu)$ given by
$(x_g)_g \sim (y_g)_g$ iff $\exists F \subset G$ finite with
$x_g=y_g, \forall g\in G \setminus F$, $\Pi_{g\in F} \mu_0(x_g) =
\Pi_{g\in F} \mu_0(y_g)$. Denote $R_0 = M(\Cal R_0)$ the von
Neumann algebra associated with $\Cal R_0$, isomorphic to the
hyperfinite II$_1$ factor. Note that the Bernoulli shift action
$\sigma_0$ of $G$ on $(X,\mu)$ leaves $\Cal R_0$ invariant, thus
implementing an action $\sigma$ of $G$ on $R_0$ and an equivalence
relation $\Cal S=\Cal R_0 \rtimes G$ on $(X,\mu)$. If we let $\Cal
R = \Cal R_0 \rtimes \Bbb Z^2$, $R=R_0\rtimes_{\sigma_0} \Bbb Z^2$
and $\sigma'$ be the action of $SL(2,\Bbb Z)$ on $R$ implemented
by $\sigma$, then the action $\sigma_{n,1}$ in the proof of 5.4
and 5.5 coincides with the restriction of $\sigma'$ to  a subgroup
$\Bbb F_n \subset SL(2,\Bbb Z)$.

$2^\circ$. In ([Fu]), A. Furman provided examples of ergodic
countable measure preserving equivalence relations $\Cal R$ with
the property that $\Cal R^t$ is not implementable by a free
ergodic action of a countable group, $\forall t > 0$, thus solving
a longstanding problem in ([FM]). It seems that the equivalence
relations constructed in part 1$\circ$ above and in the proof of
5.4 give another class of examples of equivalence relations with
this property, whenever $S \neq \{1\}$. However, we were not able
to prove this.

\heading 6. The factors $N \rtimes_\sigma G$  and the class HT
\endheading

In this Section we relate the class of cross product factors $N
\rtimes_\sigma G$ studied in this paper with the HT factors in
([Po3]), proving that the two classes are essentially disjoint.
More precisely, we prove that if $\sigma$ is sub-malleable mixing,
then $N \rtimes_\sigma G$ cannot be HT$_{_{ws}}$, a property in
between the HT and HT$_{_{s}}$ conditions considered in ([Po3]).
In the process, we obtain  two new examples of free ergodic
measure preserving actions of $\Bbb F_n$ on the probability space,
$2\leq n \leq \infty$, non orbit equivalent to the three actions
constructed in (8.7 of [Po3]).

In fact, one can show that these two actions are not orbit
equivalent to the ``residually finite'' actions of $\Bbb F_n$
constructed in ([Hj]) either, thus providing the 5'th and 6'th
distinct actions of $\Bbb F_n$. We will discuss this fact in a
forthcoming paper. In this respect, we should mention that by
([GaPo]) one knows that there exist uncountably many non orbit
equivalent free ergodic measure preserving actions of $\Bbb F_n$
on the probability space, for each $2 \leq n \leq \infty$.
However, this is an ``existence'' result, so it is still quite
interesting to prove that certain specific actions of $\Bbb F_n$
are not orbit equivalent.

\vskip .05in \noindent {\it 6.1 Definition}. A Cartan subalgebra
$A$ in a II$_1$ factor $M$ is called a {\it HT$_{_{ws}}$ Cartan
subalgebra} if $M$ has the property H relative to $A$ and there
exists a von Neumann subalgebra $A_0 \subset A$ such that:
$A_0'\cap M = A$, $\Cal N_M(A_0)''=M$, $A_0 \subset M$ is a rigid
inclusion (see [Po3] for the definitions of relative property H
and rigid inclusions).

\vskip .05in
Recalling from ([Po3]) that a Cartan subalgebra $A
\subset M$ is HT$_{_{s}}$ (resp. HT) if $M$ has the property H
relative to $A$ and $A \subset M$ is rigid (resp. $\exists$ $A_0
\subset A$ with $A_0'\cap M = A$ and $A_0 \subset M$ rigid), it
follows that $A \subset M$ is HT$_{_{s}}$ implies $A$ is
HT$_{_{ws}}$, which in turn implies $A$ is HT. Note also that all
concrete examples of HT Cartan subalgebras constructed in ([Po3])
are in fact HT$_{_{ws}}$ (see the proof of Corollary 6.3 below).
In particular, this is the case with the examples coming from the
actions $\sigma_i, i=1,2,3,$ in (5.3.3$^\circ$ and 8.7 of [Po3]).

\proclaim{6.2. Theorem} Let $\Gamma$ be an infinite group with
Haagerup's compact approximation property.

$1^\circ$. If $\sigma$ is a sub-malleable mixing action of
$\Gamma$ on $(N, \tau)$ then $M=N \rtimes_\sigma \Gamma$ contains
no diffuse subalgebra $Q \subset M$ with $P=\Cal N_M(Q)''$ a
factor and $Q \subset P$ rigid. In particular, $M$ is not a
${\text{\rm HT}}_{_{ws}}$ factor.

$2^\circ$. If $\sigma : G_0 \rightarrow {\text{\rm Aut}}(L^\infty(X, \mu))$
is a free ergodic action such that $A \subset A\rtimes_\sigma G_0$
is a rigid inclusion, then for any ergodic action
$\sigma' : G_0 \rightarrow {\text{\rm Aut}}(L^\infty(X', \mu'))$,
the diagonal product action $\sigma\otimes \sigma'
: G_0 \rightarrow {\text{\rm Aut}}(L^\infty(X \times X', \mu\times \mu')$
gives rise to a $\text{\rm HT}_{_{ws}}$ Cartan subalgebra
$L^\infty(X \times X', \mu\times \mu')=A
\subset M = A
\rtimes_{\sigma\times \sigma'} G_0$.
\endproclaim
\noindent {\it Proof}. 1$^\circ$. By Theorem 4.4, if $Q \subset M$
is diffuse with $\Cal N_M(Q)''=P$ a factor and $Q \subset P$
rigid, then there exists $\theta_\beta \in$
Aut$_\beta(M;\{\tilde{\sigma}_n\}_n)$ with $\theta_\beta(P)\subset
L(G)^\beta $, for some $\beta \in S(\{\tilde{\sigma}\}_n)$,
$\{\tilde{\sigma}\}_n$ being a given sequence of gauged extensions
for $\sigma$.

By (4.5.2$^\circ$ in [Po3]),
it follows that $\theta_\beta(Q) \subset L(G)^\beta$ is a rigid inclusion.
By cutting with an appropriate  projection $q\in Q$
of sufficiently small trace, it follows that we have  a
rigid inclusion $Q_0 \subset pL(G)p$,
where $1_{Q_0}=p = \theta_\beta(q) \in L(G)$ and $Q_0 =\theta_\beta(qQq)$.

On the other hand, by (3.1 in [Po3]) $M$ has the property H
relative to $N$, so there exists a sequence of unital
$N$-bimodular completely positive maps $\Phi_n$ on $M$ such that
$\underset n \rightarrow \infty \to \lim \|\Phi_n(x)-x\|_2=0$,
$\forall x\in M$, $\tau \circ \Phi_n = \tau$ and
$\Phi_{n|\ell^2(\Gamma)}$ compact, $\forall n$. By the rigidity of
$Q_0 \subset pMp$ applied to the completely positive maps
$p\Phi_n(p\cdot p)p$ on $pMp$, it follows that $\underset n
\rightarrow \infty \to \lim \|p\Phi_n(u)p-u\|_2=0$ uniformly for
$u \in \Cal U(Q_0)$. Since, we also have $\underset n \rightarrow
\infty \to \lim \|\Phi_n(p)-p\|_2=0$, it follows that there exists
$n$ such that $\Phi=\Phi_n$ satisfies
$$
\|\Phi(u)-u\|_2 \leq \tau(p)/2, \forall u\in \Cal U(Q_0) \tag 6.2.1
$$

But $Q_0$ diffuse implies there exists a unitary element $v \in
Q_0$ such that $\tau(v^m)=0, \forall m\neq 0$. Since $\Phi_n$ is
compact on $\ell^2(\Gamma)=L^2(L(\Gamma), \tau)$ and $v^m$ tends
weakly to $0$, it follows that $\underset m \rightarrow \infty \to
\lim \|\Phi(v^m)\|_2=0$, and for $m$ large enough, $u=v^m$
contradicts $6.2.1$.

$2^\circ$. If we put $A_0 = L^\infty(X, \mu)\subset A$ then by
construction we have $A_0 \subset M$ rigid inclusion and $A_0'\cap
M = A$. Since $\Cal N_M(A_0)$ contains both $\Cal U(A)$ and the
canonical unitaries $\{u_g\}_g$ implementing the cross product, we
also have $\Cal N_M(A_0)'' = M$. \hfill Q.E.D.

\proclaim{6.3. Corollary} Let $\Gamma$ be an ICC group with
Haagerup's compact approximation property. Assume $\Gamma$ can act
outerly and ergodically on an infinite abelian group $H$ such that
the pair $(H \rtimes \Gamma, H)$ has the relative property $(T)$.
Denote by $\sigma_1$ the action of $\Gamma$ on $L^\infty (X_1,
\mu_1) \simeq L(H)$ implemented by the action of $\Gamma$ on $H$.
Let $\sigma_0$ be a (classic) Bernoulli shift action of $\Gamma$
and $\sigma'$ an ergodic but not strongly ergodic action of
$\Gamma$ on $L^\infty(X', \mu')$ (cf. $\text{\rm [CW]}$). Let
$\sigma_2, \sigma_3$ denote the product actions of $\Gamma$ given
by $\sigma_2(g) = \sigma_1(g) \otimes \sigma_0(g)$, $\sigma_3(g) =
\sigma_1(g) \otimes \sigma'(g)$. Then $\sigma_i, 0\leq i \leq 3,$
are mutually non orbit equivalent.

Moreover, if $\Gamma$ has an infinite amenable quotient $\Gamma'$,
with $\alpha: \Gamma \rightarrow \Gamma'$ the corresponding
quotient map, and if $\sigma_0'$ denotes a Bernoulli shift action
of $\Gamma'$ on the standard probability space, then we can take
$\sigma'$ above of the form $\sigma' = \sigma'_0 \circ \alpha$ and
the product action $\sigma_4 = \sigma_0 \otimes (\alpha \circ
\sigma_0')$ is not orbit equivalent to $\sigma_i, 0\leq i \leq 3$.
Also, $\sigma_1$ has at most countable ${\text{\rm Out}}$-group
while $\sigma_0, \sigma_2, \sigma_3, \sigma_4$ have uncountable
${\text{\rm Out}}$-group.
\endproclaim
\noindent {\it Proof}. Note first that the condition on $H \subset
H \rtimes \Gamma$ in the hypothesis of the statement is equivalent
to $L(H \rtimes \Gamma)$ being a $\text{\rm HT}_{_{s}}$ factor
with $L(H)$ its $\text{\rm HT}_{_{s}}$ Cartan subalgebra, in the
sense of $\text{\rm [Po3]}$.

By 1.6.1, both $\sigma_0$ and $\sigma_4$ implement malleable
mixing integral preserving actions on the standard non-atomic
abelian von Neumann algebra $A$. By 6.2.2$^\circ$, $A
\rtimes_{\sigma_i} \Gamma, i=1,2,3,$ are ${\text{\rm HT}}_{_{ws}}$
factors, while by 6.2.1$^\circ$ it follows that $A
\rtimes_{\sigma_j} \Gamma$ are not isomorphic to $A
\rtimes_{\sigma_i} \Gamma$ for $j=0,1, i=1,2,3$. Thus, by ([FM]),
$\sigma_j$  are not orbit equivalent to $\sigma_i, $ for $j=0,1$
and $i=1,2,3.$ Also, since $\Gamma'$ is amenable, the action
$\sigma_0'$ is not strongly ergodic, so $\sigma_4$ is not strongly
ergodic either. Thus $\sigma_0$ and $\sigma_4$ are not orbit
equivalent. \hfill Q.E.D.

\proclaim{6.4. Corollary} For each $2\leq n \leq \infty$,
$\{\sigma_i\}_{0\leq i\leq 4}$ give five free ergodic measure
preserving actions of $\Bbb F_n$ on the standard probability space
that are mutually non orbit equivalent. Moreover, $\sigma_1$ has
at most countable ${\text{\rm Out}}$-group while $\sigma_0,
\sigma_2, \sigma_3, \sigma_4$ have uncountable ${\text{\rm
Out}}$-group.
\endproclaim
\noindent
{\it Proof}. This is immediate from 6.3, since
$\Bbb F_n$ does have infinite amenable quotients.
\hfill Q.E.D.

\head References\endhead

\item{[AW]} H. Araki, J. Woods: {\it A classification of factors},
Publ. Res. Math. Sci., Kyoto Univ., {\bf 6} (1968), 51-130.

\item{[Bu]} M. Burger, {\it Kazhdan constants for} $SL(3,\Bbb Z)$,
J. reine angew. Math., {\bf 413} (1991), 36-67.

\item{[Ch]}
E. Christensen: {\it Subalgebras of a finite algebra},
Math. Ann. {\bf 243} (1979), 17-29.

\item{[Cho]} M. Choda: {\it Group factors of the Haagerup type},
Proc. Japan Acad., {\bf 59} (1983), 174-177.

\item{[C1]} A. Connes: {\it A type II$_1$
factor with countable fundamental group}, J. Operator
Theory {\bf 4} (1980), 151-153.

\item{[C2]} A. Connes: {\it Une classification des facteurs de type III},
Ann. Ec. Norm. Sup. {\bf 6} (1973), 133-252.

\item{[C3]} A. Connes: {\it Almost periodic states and factors of type
III$_1$},
J. Funct. Anal. {\bf 16} (1974), 415-455.

\item{[CJ]} A. Connes, V.F.R. Jones: {\it Property T
for von Neumann algebras}, Bull. London Math. Soc. {\bf 17} (1985),
57-62.

\item{[CS]} A. Connes, E. St\o rmer: {\it Entropy for automorphisms
of II$_1$ von Neumann
algebras}, Acta Math. {\bf 134} (1974), 289-306.

\item{[DeKi]} C. Delaroche, Kirilov: {\it Sur les relations entre
l'espace dual d'un groupe et la structure de ses sous-groupes fermes},
Se. Bourbaki, 20'eme ann\'ee, 1967-1968, no. 343, juin 1968.

\item{[Di]} J. Dixmier: ``Les alg\`ebres d'op\'erateurs sur l'espace
Hilbertien
(Alg\`ebres de von Neumann)'', Gauthier-Villars, Paris, 1957, 1969.

\item{[FM]} J. Feldman, C.C. Moore: {\it Ergodic equivalence relations,
cohomology, and von Neumann algebras I, II}, Trans. Amer. Math.
Soc. {\bf 234} (1977), 289-324, 325-359.

\item{[Fu]} A. Furman: {\it Orbit equivalence
rigidity}, Ann. of Math. {\bf 150} (1999), 1083-1108.

\item{[Ga]} D. Gaboriau: {\it Invariants $\ell^2$ de r\'elations
d'\'equivalence et de groupes}, Publ. Math. I.H.\'E.S. {\bf 95}
(2002), 93-150.

\item{[GaPo]} D. Gaboriau, S. Popa: {\it An Uncountable Family of
Non Orbit Equivalent Actions of $\Bbb F_n$}, preprint,
math.GR/0306011, to appear in J. of Amer. Math. Soc.

\item{[GeGo]} S. Gefter, V.  Golodets: {\it Fundamental groups for
ergodic actions and actions with unit fundamental groups}, Publ.
RIMS, Kyoto Univ. {\bf 24} (1988), 821-847.

\item{[GoNe]} V. Y. Golodets, N. I. Nesonov: T{\it -property and
nonisomorphic factors of type} II {\it and} III, J. Funct. Analysis
{\bf 70} (1987), 80-89.

\item{[H]} U. Haagerup: {\it An example of non-nuclear C$^*$-algebra
which has the metric approximation property}, Invent. Math.
{\bf 50} (1979), 279-293.

\item{[dHV]} P. de la Harpe, A. Valette: ``La propri\'et\'e T
de Kazhdan pour les
groupes localement compacts'', Ast\'erisque {\bf 175} (1989).

\item{[Hj]} G. Hjort: {\it A converse to Dye's Theorem},
UCLA preprint, September 2002.

\item{[Jo]} P. Jolissaint: {\it On the relative property T},
preprint 2001.

\item{[J1]} V.F.R. Jones : {\it A converse to Ocneanu's theorem},
J. Operator Theory {\bf 4} (1982), 21-23.

\item{[J2]} V.F.R. Jones : {\it Index for subfactors}, Invent. Math.
{\bf 72} (1983), 1-25.

\item{[J3]} V.F.R. Jones : {\it Ten problems}, in ``Mathematics:
perspectives and frontieres'', pp. 79-91, AMS 2000, V. Arnold, M. Atiyah,
P. Lax, B. Mazur Editors.

\item{[JPo]} V.F.R. Jones, S. Popa : {\it Some properties
of MASA's in factors}, in
``Invariant subspaces and other topics'', pp. 89-102,
Operator Theory: Adv. Appl., by Birkhauser, Boston, 1982.

\item{[K]} R.V. Kadison: {\it Problems on von Neumann algebras},
Baton Rouge Conference 1967, unpublished.

\item{[Ka]} D. Kazhdan: {\it Connection of the dual space of a group
with the structure of its closed subgroups}, Funct. Anal. and its Appl.,
{\bf1} (1967), 63-65.

\item{[Ma]} G. Margulis: {\it Finitely-additive invariant measures
on Euclidian spaces}, Ergodic. Th. and Dynam. Sys. {\bf 2} (1982),
383-396.

\item{[MS]} N. Monod, Y. Shalom:
{\it Orbit equivalence rigidity and bounded cohomology},
Preprint 2002.

\item{[MvN1]} F. Murray, J. von Neumann:
{\it Rings of operators}, Ann. Math. {\bf 37}
(1936), 116-229.

\item{[MvN2]} F. Murray, J. von Neumann:
{\it On rings of operators IV}, Ann. Math. {\bf 44} (1943),
716-808.

\item{[Po1]} S. Popa: {\it Some rigidity results for non-commutative
Bernoulli shifts}, MSRI preprint 2001-005).

\item{[Po2]} S. Popa: {\it Correspondences},
INCREST preprint 1986, unpublished.

\item{[Po3]} S. Popa: {\it On a class of type II$_1$ factors with
Betti numbers invariants}, MSRI preprint no 2001-024, revised
math.OA/0209130, to appear in Ann. of Math.

\item{[Po4]} S. Popa: {\it On the fundametal group of type} II$_1$
{\it factors}, Proc. Nat. Acad. Sci. {\bf 101} (2004),
723-726. (math.OA/0210467)

\item{[Po5]} S. Popa: {\it
Free independent sequences in type} $II_1$ {\it factors and related
problems}, Asterisque {\bf 232} (1995), 187-202.

\item{[Po6]} S. Popa: {\it Strong rigidity of} II$_1$ {\it factors
arising from malleable actions of w-rigid groups} II,
(math.OA/0407137) and III (in preparation).

\item{[Po7]} S. Popa: {\it Strong rigidity of} II$_1$ {\it factors
arising from malleable actions of w-rigid groups} III, in preparation.

\item{[P]} R. Powers: {\it Representation of uniformly hyperfinite
algebras and their associated von Neumann rings}, Ann. Math.
{\bf 86} (1967), 138-171.

\item{[PSt]} R. Powers, E. St\o rmer: {\it Free
states of the canonical anticommutation relations},
Comm. Math. Phys. {\bf 16} (1970), 1-33.

\item{[Sa]} Sakai: ``C$^*$-algebras and W$^*$-algebras'', Springer-Verlag,
Berlin-Heidelberg-New York, 1971.

\item{[T1]} M. Takesaki: {\it The structure of a von Neumann
algebra with a homogeneous periodic state},
Acta Math. {\bf 131} (1973), 281-308.

\item{[T2]} M. Takesaki: {\it Conditional expectation
in von Neumann algebra}, J. Fnal. Analysis {\bf 9} (1972), 306-321.

\item{[T3]} M. Takesaki: ``Theory of Operator Algebras II'',
Encyclopedia of Mathematical Sciences {\bf 125}, Springer-Verlag,
Berlin-Heidelberg-New York, 2002.

\item{[Va]} A. Valette: {\it Group pairs with relative property} (T)
{\it from arithmetic lattices}, preprint 2004 (preliminary version
2001).

\item{[Zi]} R. Zimmer: ``Ergodic theory and semisimple groups'',
Birkha\"user-Verlag, Boston 1984.

\enddocument